\documentclass[letterpaper,9pt]{IEEEtran}

\usepackage{amsmath,amscd,amsbsy,amssymb,latexsym,url,bm,amsthm}
\usepackage{array}
\usepackage{algorithm}
\usepackage{algorithmic}
\usepackage{booktabs}
\usepackage{cases}
\usepackage{cite}
\usepackage{etoolbox}
\preto\subequations{\ifhmode\unskip\fi} 
\usepackage[table]{xcolor}
\usepackage{graphicx}
\usepackage{hyperref}
\hypersetup{hidelinks}
\usepackage{multirow,makecell}
\usepackage{mathtools}
\usepackage{subeqnarray}
\usepackage[caption=false,font=footnotesize,subrefformat=parens,labelformat=parens]{subfig}
\usepackage{tabularx}
\usepackage{threeparttable}
\usepackage{tikz}
\usepackage{balance}
\usepackage[capitalize]{cleveref}
\allowdisplaybreaks

\newtheorem{assumption}{Assumption}
\newtheorem{theorem}{Theorem}
\newtheorem{lemma}{Lemma}

\newtheorem{remark}{Remark}

\crefname{assumption}{Assumption}{Assumptions}
\crefname{corollary}{Corollary}{Corollaries}
\crefname{theorem}{Theorem}{Theorems}
\crefname{lemma}{Lemma}{Lemmas}
\crefalias{lemma}{lemma}

\newcommand{\bigO}[1]{\ensuremath{\mathop{}\mathopen{}\mathcal{O}\mathopen{}\left(#1\right)}}

\renewcommand{\bar}[1]{\mkern 1.5mu\overline{\mkern-1.5mu#1\mkern-1.5mu}\mkern 1.5mu}
\newcommand{\ud}{\mathrm{d}}
\newcommand{\RNum}[1]{\uppercase\expandafter{\romannumeral #1\relax}}
\newcommand{\numcircled}[1]{\tikz[baseline=(char.base)]{
            \node[shape=circle,draw,inner sep=0.4pt] (char) {#1};}}
\newcommand{\obj}{\tilde{\Phi}}

\newcommand{\rev}[1]{{\color{black} #1}} 
\newcommand{\revise}[1]{{\color{black} #1}} 

\DeclareMathOperator*{\argmin}{arg\,min}
\DeclareMathOperator{\diag}{diag}
\DeclareMathOperator{\E}{\mathbb{E}}
\DeclareMathOperator{\proj}{Proj_{\mathcal{U}}}
\DeclareMathOperator{\reg}{Reg_\mathnormal{T}^\mathnormal{d}}

\newcolumntype{Y}{>{\centering\arraybackslash}X}

\makeatletter
\def\raisedotfill{%
  \leavevmode
  \cleaders \hb@xt@ .44em{\hss\raise0.5ex\hbox{.}\hss}\hfill
  \kern\z@}
\pretocmd\@bibitem{\color{black}\csname keycolor#1\endcsname}{}{\fail}
\newcommand\citecolor[1]{\@namedef{keycolor#1}{\color{blue}}}
\let\c@lemma\c@theorem
\makeatother

\definecolor{lightgreen}{RGB}{225, 245, 230}

\newskip\origabovedisplayskip
\newskip\origbelowdisplayskip
\newskip\origabovedisplayshortskip
\newskip\origbelowdisplayshortskip

\AtBeginDocument{%
  \origabovedisplayskip=\abovedisplayskip
  \origbelowdisplayskip=\belowdisplayskip
  \origabovedisplayshortskip=\abovedisplayshortskip
  \origbelowdisplayshortskip=\belowdisplayshortskip

  \setlength{\abovedisplayskip}{5pt plus 2pt minus 2.5pt}
  \setlength{\belowdisplayskip}{5pt plus 2pt minus 2.5pt}
  \setlength{\abovedisplayshortskip}{3pt plus 2pt minus 2pt}
  \setlength{\belowdisplayshortskip}{3pt plus 2pt minus 2pt}
}

\newcommand{\restoreIEEEdisplayskips}{%
  \abovedisplayskip=\origabovedisplayskip
  \belowdisplayskip=\origbelowdisplayskip
  \abovedisplayshortskip=\origabovedisplayshortskip
  \belowdisplayshortskip=\origbelowdisplayshortskip
}

\newif\ifshowtracking
\showtrackingfalse

\newif\ifarxiv
\arxivtrue

\begin{document}
    \title{Gray-Box Nonlinear Feedback Optimization}
    \author{Zhiyu He, Saverio Bolognani, Michael Muehlebach, and Florian D{\"o}rfler
        \thanks{This work was supported by the Max Planck ETH Center for Learning Systems, the German Research Foundation, and the Swiss National Science Foundation through NCCR Automation under Grant 51NF40\_225155.}
        \thanks{
        Zhiyu He, Saverio Bolognani, and Florian D{\"o}rfler are with the Automatic Control Laboratory, ETH Z{\"u}rich, 8092 Z{\"u}rich, Switzerland (email: \{zhiyhe, bsaverio, dorfler\}@ethz.ch). Michael Muehlebach is with the Max Planck Institute for Intelligent Systems, 72076 T{\"u}bingen, Germany (email: michael.muehlebach@tuebingen.mpg.de). 
            }
            }
    \maketitle

    \begin{abstract}
        Feedback optimization enables autonomous optimality seeking of a dynamical system through its closed-loop interconnection with iterative optimization algorithms. Among various iteration structures, model-based approaches require the input-output sensitivity matrix of the system to construct gradients, whereas model-free approaches eliminate this need by estimating gradients from real-time objective evaluations. These approaches offer complementary benefits in sample efficiency and accuracy against model mismatch, i.e., sensitivity errors. To achieve balanced closed-loop performance, we propose a gray-box feedback optimization controller, featuring systematic incorporation of approximate sensitivities into model-free updates via a tunable convex combination. We provide unified performance characterizations covering different approaches. We elucidate how cumulative sensitivity errors (model-based) and variances due to stochastic exploration (model-free) shape the closed-loop behavior and induce a trade-off between iteration and dimensional dependence. The proposed controller retains sample efficiency and provable (local) optimality for nonconvex problems despite inaccurate sensitivities. We further develop and characterize a running gray-box controller that handles constrained time-varying problems with changing objectives and steady-state input-output maps.
    \end{abstract}


\section{Introduction}\label{sec:introduction}
Efficient steady-state operation is crucial for engineering systems, e.g., power grids, process control systems, and communication networks \cite{simonetto2020time}. To this end, numerical optimization provides tractable solutions based on an explicit problem formulation, e.g., an economic efficiency objective together with constraints that encode the input-output map and forecasts or anticipated statistics of disturbances. While these solutions are ready to be implemented in an offline and feedforward fashion, their effectiveness may be jeopardized if the complex plant and unmeasured disturbances are not accurately modeled. Given the extensive use of feedback control for adaptation and robustness, it is desirable to pursue an optimal steady-state behavior through a closed-loop paradigm.

\subsection{Related Work}
Closed-loop optimization has been explored from multiple perspectives. A core principle is to interact with plant dynamics, learn from feedback, and improve control strategies. Typical examples focusing on cumulative performance include model predictive control and reinforcement learning. Iterative adjustments of strategies based on feedback are pivotal in handling inaccurate models or unknown dynamics and navigating through vast policy spaces. 

In terms of closed-loop steady-state optimization, extremum seeking is an effective solution free of any model information. The central idea is to add random exploration signals (e.g., sinusoids) and use averaging to construct an update direction. Extremum seeking requires exploration signals of well-selected frequencies to avoid mutual interference, and therefore it is typically applied to low-dimensional problems. Constraint satisfaction can be achieved through various techniques, e.g., penalty or barrier functions, saddle-point dynamics, optimization on manifolds, and projection maps \cite{scheinker2024100}.

Feedback optimization\cite{simonetto2020time} emerges as a promising paradigm for steady-state optimization of dynamical systems. Compared to extremum seeking, this paradigm systematically handles high-dimensional objectives and constraints. The key insight is to implement optimization-based iterations as feedback controllers, thereby driving a stable plant to an optimal steady state. Thanks to the use of real-time measurements, feedback optimization bypasses the need to explicitly access the complex input-output map and the unknown exogenous disturbance. Moreover, it handles (partially) unknown objectives by using estimation via parametric or nonparametric models based on noisy function evaluations\cite{simonetto2021personalized,cothren2022online}. Provided that the controller gain is low (i.e., a time-scale separation holds)\cite{simpson2021analysis}, closed-loop stability, optimality, and constraint satisfaction are guaranteed\cite{bianchin2021time}. This closed-loop structure also facilitates tracking the trajectory of time-varying optimal solutions in non-stationary environments\cite{belgioioso2022online,cothren2022online} and proves effective in industry applications\cite{simonetto2020time}.


The manifold benefits of feedback optimization rely on the premise that the steady-state input-output sensitivity matrix of the plant is available. This premise stems from using the chain rule to formulate the gradient-based update direction of the controller. 
In practice, the challenges of modeling complex, large-scale, and poorly known systems can render accurate sensitivities elusive, thereby causing closed-loop sub-optimality, constraint violation, or instability. Two streams of strategies have been explored to address this issue.

One stream is \emph{model-based}, since the key model information (i.e., sensitivity) is learned from offline data or online interactions.
Behavioral systems theory enables data-driven representations of the sensitivity of a linear time-invariant plant via its historical trajectories\cite{bianchin2021online,nonhoff2022online}.
Moreover, recursive least-squares estimation allows learning sensitivities in an online fashion based on streaming data of inputs and outputs\cite{picallo2021adaptive,dominguez2023online}. 
Nonetheless, if the sensitivity is not learned fast and accurately enough, the closed-loop performance may experience considerable sub-optimality, see \cite[Section~VI]{he2022model}.

Another stream to address unknown sensitivities is \emph{model-free}, in that learning sensitivities is avoided altogether. 
This goal is achieved by employing iterative optimization schemes without gradient evaluations. 
Some methods embrace the probabilistic framework of Bayesian optimization and adjust the input as the maximizer of an acquisition function \cite{krishnamoorthy2023model}.
A different strategy better suited for high-dimensional objectives is to leverage zeroth-order optimization\cite{nesterov2017random} and construct stochastic gradient estimates from function evaluations \cite{poveda2017robust,chen2020model,tang2023zeroth,chen2025continuous,he2022model}. 
\revise{Overall, the stochasticity of gradient estimates renders the corresponding convergence guarantee dependent on the problem dimension, resulting in an increased actuation count relative to model-based approaches.}

\subsection{Motivation}\label{subsec:motivations}
Models encode useful structural information, contributing to the high sample efficiency (i.e., fast convergence) of model-based controllers. This advantage, however, comes with stringent accuracy requirements for models that must be formulated or learned. Model-free operations are attractive because they offer provable guarantees without resorting to complex models. Nevertheless, they can be less sample-efficient due to stochastic exploration or be restricted to certain classes of problems (e.g., in low-dimensional spaces).

Given such complementary benefits, it is promising to develop \emph{gray-box} approaches to achieve the best of both worlds, as evidenced by examples in reinforcement learning, online control\cite{qu2021exploiting}, predictive control, and stabilization\cite{li2023certifying}. 
Some methods are built upon model-based pipelines and introduce model-free blocks for inference or improvement, including estimating initial states, generating feedforward inputs\cite{ma2023reinforcement}, and learning terminal costs and constraints. \revise{Others augment model-free pipelines with model-based priors in a sequential manner \cite{qu2021exploiting} or utilize synthetic data generated from transition models to enhance model-free training.} 
Furthermore, learning-augmented control\cite{li2023certifying} combines a model-based (albeit sub-optimal) policy with a machine-learned policy, a viewpoint reminiscent of multi-model adaptive control. 
However, the tuning of the combination coefficients therein requires accessing system matrices or errors in state estimation, which may not always be available in applications.


\revise{Model-based pipelines are preferable given highly accurate models, whereas model-free pipelines may have an edge in the absence of model information and structure. For steady-state performance optimization, we focus on an intermediate \emph{gray-box} regime, where approximate input-output sensitivity matrices of a physical plant are useful but not highly accurate. In practice, such sensitivities are obtained through prior knowledge, first-principles models, or recursive estimation\cite{ma2023reinforcement}. We leverage approximate sensitivities to design gray-box feedback optimization controllers, featuring principled interpolation, unified analysis, and explicit performance trade-offs.}

\subsection{Contributions}
\revise{
We develop gray-box feedback optimization controllers that utilize approximate sensitivities for closed-loop steady-state performance optimization.
Our contributions are summarized as follows.


\begin{itemize}
	\item We propose gray-box controllers that drive a stable nonlinear system to optimal steady-state operating points. When implemented in closed loop, these controllers leverage real-time output measurements and iteratively adjust inputs by interpolating between model-based inexact gradients (constructed from approximate sensitivities) and model-free gradient estimates through a convex combination scheme.

	\item Owing to the above principled interpolation, our gray-box controllers encompass model-based, model-free, and sequential (i.e., first model-based, and then model-free, or vice versa) pipelines as special cases, thereby providing unified formulations for the design and analysis of feedback optimization. Furthermore, the proposed controllers are flexible in handling approximate sensitivity matrices of varying qualities, by virtue of an interpolation that carefully modulates emphasis between model-based and model-free directions over time.


	\item We establish closed-loop performance guarantees for static nonconvex problems and time-varying convex problems. To this end, we decouple the effects of sensitivity errors (model-based) and stochastic exploration (model-free) on the corresponding performance measures, i.e., the average second moment of gradients and the dynamic regret. The proposed controllers achieve balanced closed-loop performance in sample efficiency (with trade-offs in iteration and dimensional dependence) and optimality despite (persistent) sensitivity errors.
\end{itemize}
}

The rest of this article is organized as follows. In \cref{sec:formulation}, we introduce the problem of interest and some preliminaries. \cref{sec:design} presents the design of our gray-box controller. The performance guarantees in a static and unconstrained setting are established in \cref{sec:analysis}. \cref{sec:tracking} explores the extension to handle time-varying problems with input constraints. We perform numerical evaluations in \cref{sec:experiment}. Finally, \cref{sec:conclusion} concludes this article. 
\ifarxiv
All proofs are deferred to the appendix.
\else
All proofs are deferred to the online report \cite{he2024gray}.
\fi

\section{Problem Formulation and Preliminaries}\label{sec:formulation}
\subsection{Problem Formulation}\label{subsec:formulation}


We aim to find an input $u$ to optimize the steady-state operation of a nonlinear physical plant, i.e.,
\begin{equation}\label{eq:opt_original}
\begin{split}
    \min_{u\in \mathbb{R}^p, y\in \mathbb{R}^q} \quad & \Phi(u,y) \\ 
    \text{s.t.} \quad & y = h(u,d).
\end{split}
\end{equation}
In \eqref{eq:opt_original}, the objective $\Phi: \mathbb{R}^p \times \mathbb{R}^q \to \mathbb{R}$ is a nonconvex function of the input $u \in \mathbb{R}^p$ and the steady-state output $y \in \mathbb{R}^q$ of a fast stable nonlinear plant. Further, $d \in \mathbb{R}^r$ is a fixed unknown exogenous disturbance. The constraint of \eqref{eq:opt_original} encodes the steady-state input-output map $h:\mathbb{R}^{p} \times \mathbb{R}^{r} \to \mathbb{R}^q$ of this plant, i.e.,
\begin{equation}\label{eq:sys_map}
    y = h(u,d).
\end{equation}
While in \eqref{eq:opt_original} we consider a plant abstracted by an algebraic map to streamline the presentation of guarantees, the online nature of feedback optimization controllers allows adaptation to handle stable nonlinear dynamics. We will outline theoretical extensions in \cref{sec:analysis} and present numerical evaluations in \cref{sec:experiment}. The broader scenario involving input constraints and time-varying disturbances will be examined in \cref{sec:tracking}. Our assumptions are as follows.

\begin{assumption}\label{assump:sys_map}
    The steady-state map $h(u,d)$ is differentiable with respect to $u$ \rev{for any given $d \in \mathbb{R}^r$.}
\end{assumption}

\begin{assumption}\label{assump:objective}
    The reduced objective $\obj(u) \triangleq \Phi(u, h(u,d))$ is $M$-Lipschitz and $L$-smooth (i.e., with $L$-Lipschitz gradients) \revise{in $u$ for any given $d \in \mathbb{R}^r$}. It satisfies $\obj^* \triangleq \inf_{u\in \mathbb{R}^p} \obj(u) > -\infty$.
\end{assumption}

\begin{assumption}\label{assump:objective_Lipschitz}
    The function $\Phi(u,y)$ is $M_\Phi$-Lipschitz in $y$.
\end{assumption}

\cref{assump:sys_map} is typical and encompasses a broad class of systems. For instance, a stable linear dynamical system admits a linear differentiable steady-state map\cite{simpson2021analysis}. For a nonlinear dynamical system, the (local) existence and differentiability of its steady-state map can be ensured by the implicit function theorem (see \cite[Theorem~1B.1]{dontchev2009implicit}) with suitable conditions on system dynamics. The above properties of the objective function are relatively weak, commonly assumed (e.g., \cite{nesterov2017random,zhang2022new,tang2023zeroth}), and satisfied in many applications\cite{simonetto2020time}.


A tempting solution to problem~\eqref{eq:opt_original} is to directly use numerical optimization solvers. Nonetheless, solvers require the explicit steady-state map and the exact value of the disturbance. These requirements can be hard to satisfy when complex systems and unknown disturbances are present. Hence, we pursue a feedback optimization controller that utilizes real-time output measurements to iteratively drive the plant~\eqref{eq:sys_map} to an optimal operating point quantified by \eqref{eq:opt_original}.


Model-based feedback optimization controllers\cite{belgioioso2022online,picallo2021adaptive,dominguez2023online} learn and use the input-output sensitivity $\nabla_u h(u,d)$ of the plant \eqref{eq:sys_map} and iteratively update inputs by following the gradient of the reduced objective $\obj(u)$. After invoking the chain rule, their update rule reads
\begin{equation}\label{eq:model-based-FO} 
    u_{k+1} = u_k \!-\! \eta (\nabla_u \Phi(u_k,y_k) \!+\! \nabla_u h(u_k,d) \nabla_y \Phi(u_k,y_k)),
\end{equation}
where $\eta>0$ is a step size, $\nabla_u h(u_k,d)$ and $y_k$ are the sensitivity and the output of the plant \eqref{eq:sys_map} evaluated at time $k\in \mathbb{N}$, respectively, \revise{and $\nabla_u \Phi$ and $\nabla_y \Phi$ denote the partial gradients of $\Phi$ with respect to $u$ and $y$, respectively}. 
Model-free controllers\cite{chen2020model,tang2023zeroth,he2022model} bypass sensitivity information and purely rely on stochastic exploration. Their trade-offs in sample efficiency and solution accuracy are discussed in \cref{subsec:motivations}. 
In contrast, we consider a \emph{gray-box} regime, where only approximate model information, such as a sensitivity estimate derived from prior knowledge or recursive estimation, is available. \revise{We will merge approximate sensitivities into model-free updates, thus achieving balanced closed-loop performance in terms of provable accuracy and sample efficiency.}


\vspace{-2ex}
\subsection{Preliminaries of Gradient Estimation}
Various model-free feedback optimization controllers \cite{chen2020model,tang2023zeroth,he2022model} exploit zeroth-order optimization\cite{nesterov2017random}. Their underlying idea is to iteratively update in the direction of negative gradient estimates, which are constructed from function evaluations and random exploration vectors. 
Consider a function $\xi: \mathbb{R}^p \to \mathbb{R}$. Let $\delta > 0$ be a smoothing parameter. A smooth approximation $\xi_{\delta}$ of $\xi$ is 
\begin{equation}\label{eq:smooth_approx}
    \xi_{\delta}(w) = \E_{v'\sim U(\mathbb{B}_p)}[\xi(w+\delta v')],
\end{equation}
where $w\in \mathbb{R}^p$, and \rev{the expectation is taken with respect to a vector $v'$ uniformly sampled from the unit ball $\mathbb{B}_p \triangleq \{v'\in \mathbb{R}^p : \|v'\| \leq 1\}$ in $\mathbb{R}^p$}. Further, let $\mathbb{S}_{p-1} \triangleq \{v\in \mathbb{R}^p : \|v\|=1\}$ be the unit sphere in $\mathbb{R}^p$. 
The following lemma summarizes the construction of a gradient estimate and useful properties of $\xi_{\delta}$.

\begin{lemma}[\hspace{1sp}{\cite[Lemma 4.1]{gao2018information}}]\label{lem:grad_est_property}
    If $\xi: \mathbb{R}^p \to \mathbb{R}$ is $L_{\xi}$-smooth, then $\xi_{\delta}(w)$ defined in \eqref{eq:smooth_approx} is $L_{\xi}$-smooth, and
    \begin{subequations}
    \begin{align}
        \E_{v \sim U(\mathbb{S}_{p-1})}\left[\frac{p}{\delta}\xi(w+\delta v)v\right] &= \nabla \xi_{\delta}(w), \label{eq:unbiased_grad_est} \\
        |\xi_{\delta}(w) - \xi(w)| &\leq \frac{L_\xi \delta^2}{2}, \label{eq:difference_func_val} \\
        \|\nabla \xi_{\delta}(w) - \nabla \xi(w)\| &\leq \frac{L_{\xi}p\delta}{2}, \label{eq:difference_grad}
    \end{align}
    \end{subequations}
    where $w\in \mathbb{R}^p$, and $\delta > 0$. If $\xi$ is convex, then $\xi_{\delta}$ is also convex.
\end{lemma}

\rev{Lemma~\ref{lem:grad_est_property} implies that $\frac{p}{\delta}\xi(w\!+\!\delta v)v$ constructed from the function value $\xi(w\!+\!\delta v)$ and the random vector $v$ uniformly sampled from the unit sphere $\mathbb{S}_{p-1}$ is an unbiased (zeroth-order) estimate of $\nabla \xi_{\delta}(w)$.} 
\rev{The validity of \eqref{eq:unbiased_grad_est} follows from the divergence theorem, which relates the gradient of a volume integral (i.e., $\nabla \xi_{\delta}(w)$) to a surface integral (i.e., the left-hand side of \eqref{eq:unbiased_grad_est}).} Further, the closeness between the values and gradients of $\xi$ and $\xi_{\delta}$ can be adjusted to arbitrary precision via the smoothing parameter $\delta$. \rev{Finally, the randomized smoothing operation \eqref{eq:smooth_approx} preserves smoothness and convexity.} 
To connect back to feedback optimization, a gradient estimate as in \eqref{eq:unbiased_grad_est} can be employed to update the actuation $u$ in lieu of the model-based gradient as in \eqref{eq:model-based-FO}, see \cite{he2022model} for further reading.
\section{Design of the Gray-Box Controller}\label{sec:design}

\subsection{Gray-Box Feedback Optimization Controller}\label{subsec:gray_box_design}
We consider a gray-box regime, where an approximate input-output sensitivity $\hat{H}_k \in \mathbb{R}^{p \times q}$ differing from the true sensitivity $\nabla_u h(u_k,d) \triangleq H_k$ is available at time $k \in \mathbb{N}$ for a given input $u_k$. Such an approximate sensitivity can be obtained through prior knowledge, first-principles models\cite{ma2023reinforcement}, or online learning and estimation\cite{picallo2021adaptive,dominguez2023online}.

Our proposed feedback optimization controller solves \eqref{eq:opt_original} by iteratively adjusting inputs based on real-time output measurements. The update direction is constructed by adaptively fusing an inexact gradient from the approximate sensitivity $\hat{H}_k$ and a gradient estimate through stochastic exploration. The update rules are
\begin{subequations}\label{eq:hybrid_controller}
\begin{align}
    w_{k+1} &= w_k - \eta \tilde{\phi}_k, \label{eq:hybrid_controller_GD_update} \\
    \tilde{\phi}_k &= \alpha_k \tilde{\phi}_{k,1} + (1-\alpha_k) \tilde{\phi}_{k,2}, \label{eq:hybrid_controller_gradient} \\
    \tilde{\phi}_{k,1} &= \nabla_{u} \Phi(u_k,y_k) + \hat{H}_k \nabla_y \Phi(u_k,y_k), \label{eq:hybrid_controller_inexact_GD} \\
    \tilde{\phi}_{k,2} &= \frac{pv_k}{\delta} \big(\Phi(u_k,y_k) - \Phi(u_{k-1},y_{k-1})\big), \label{eq:hybrid_controller_grad_est} \\
    u_{k+1} &= w_{k+1} + \delta v_{k+1}. \label{eq:hybrid_controller_exploration}
\end{align}
\end{subequations}
In \eqref{eq:hybrid_controller}, we define $u_{-1} = y_{-1} \triangleq 0$. Further, $k \in \mathbb{N}$ is the iteration count, $w_k$ is a candidate solution, $\alpha_k\in [0,1]$ is a convex combination coefficient whose design is further specified in \cref{subsec:adaptive_weight}, $\eta > 0$ is a step size, $\delta > 0$ is a smoothing parameter, $p$ is the size of the input, and $v_0, \ldots, v_{k+1} \sim U(\mathbb{S}_{p-1})$ are independent and identically distributed (i.i.d.) random vectors sampled from the unit sphere $\mathbb{S}_{p-1}$. The role of $v_0,\ldots,v_{k+1}$ lies in exploring the objective $\obj$ around the current solutions $w_0,\ldots,w_{k+1}$, thus facilitating the construction of gradient estimates. The controller \eqref{eq:hybrid_controller} uses real-time output measurement $y_k$ as feedback from the plant \eqref{eq:sys_map}.

\begin{figure}[!tb]
    \centering
    \includegraphics[width=0.9\columnwidth]{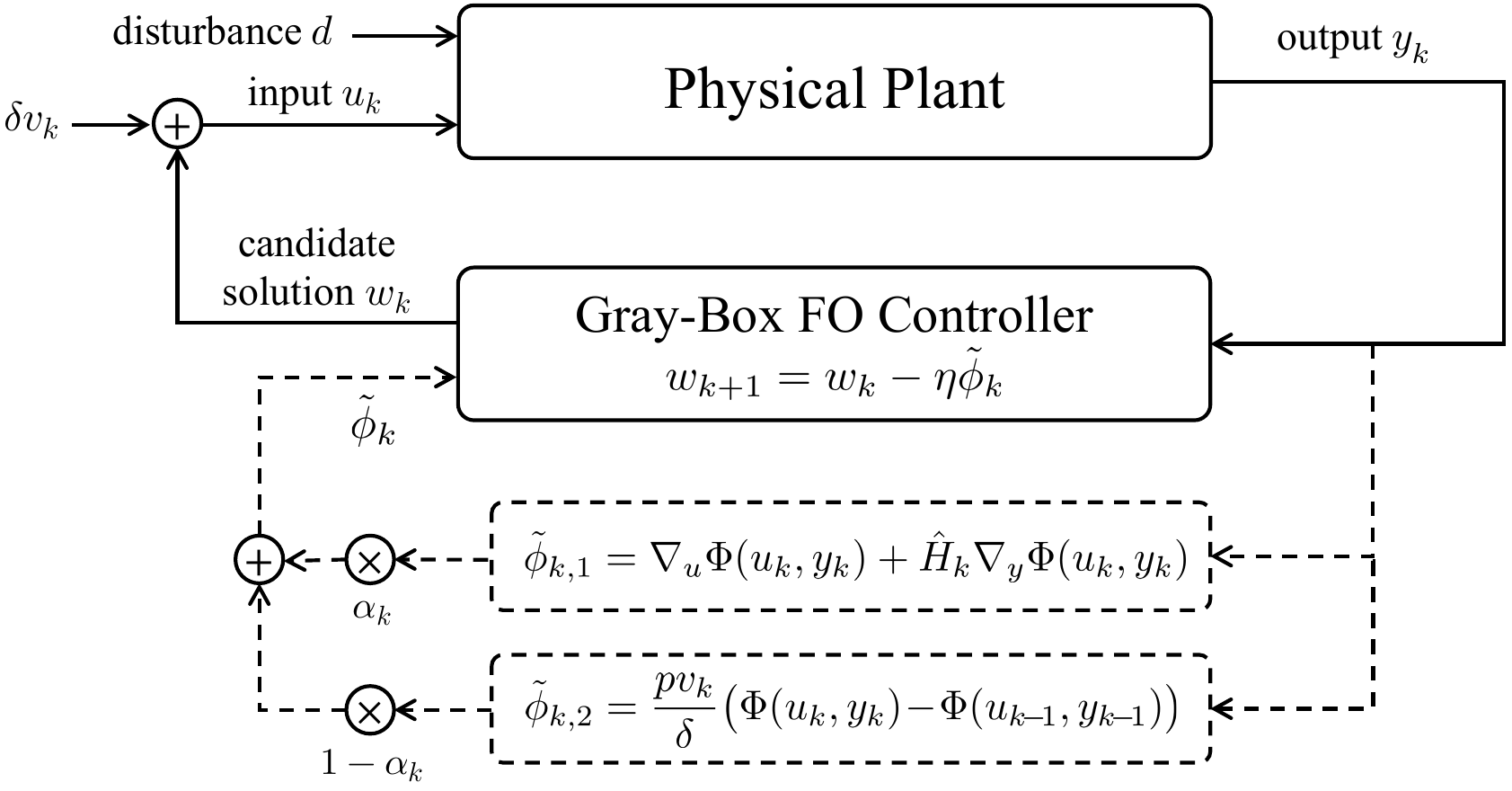}
    \vspace{-2ex}
    \caption{Interconnection of a physical plant and the gray-box feedback optimization controller \eqref{eq:hybrid_controller}.}
    \vspace{-1ex}
    \label{fig:interconnect_hybrid}
\end{figure}

The gray-box controller \eqref{eq:hybrid_controller} merges two directions via adaptive convex combination \eqref{eq:hybrid_controller_gradient}. The first (i.e., $\tilde{\phi}_{k,1}$) is an inexact gradient using $\hat{H}_k$. The second (i.e., $\tilde{\phi}_{k,2}$) is a gradient estimate constructed from the current and previous evaluations of the objective $\Phi$, which has been developed in \cite{zhang2022new} and leveraged for feedback optimization in our previous work \cite{he2022model}. \revise{Different from the estimates in \cite{chen2020model,tang2023zeroth} with two objective evaluations (i.e., two actuation steps), $\tilde{\phi}_{k,2}$ requires a single actuation per iteration. This feature aligns with online decision-making settings, where resetting and actuating the system again may be impractical.} 
Note that $\tilde{\phi}_{k,2}$ is an unbiased estimate of $\nabla \obj_\delta(w_k)$, where $\obj_\delta$ is the smooth approximation of $\obj$ as per \eqref{eq:smooth_approx}. The unbiased property follows from \eqref{eq:unbiased_grad_est}, the zero expectation of $v_k$, and the independence of $v_k$ from $v_{k-1}$ and $\Phi(u_{k-1},y_{k-1})$. 
The controller perturbs the new solution $w_{k+1}$ by an exploration noise $\delta v_{k+1}$, see \eqref{eq:hybrid_controller_exploration}, thereby exploring $\obj$ around $w_{k+1}$ to construct a gradient estimate. 
Finally, the input $u_{k+1}$ is applied to the plant \eqref{eq:sys_map}. The iterative structure also allows handling input constraints via projection, see also \cref{sec:tracking}. \cref{fig:interconnect_hybrid} illustrates the closed-loop interconnection of a physical plant and our gray-box controller \eqref{eq:hybrid_controller}.

\subsection{Adaptive Combination Coefficients}\label{subsec:adaptive_weight}
A key ingredient of our gray-box controller \eqref{eq:hybrid_controller} is the combination coefficient $\alpha_k$. This coefficient adjusts the extent of blending the approximate sensitivity $\hat{H}_k$ into model-free updates. We discuss how to tune $\alpha_k$ in different scenarios based on the quality of $\hat{H}_k$.

If we can learn the sensitivity $H_k$ sufficiently fast and accurately, e.g., via Kalman filtering\cite{picallo2021adaptive,dominguez2023online}, then model-based controllers augmented with such sensitivity estimation are favorable. Specifically, we will show in \cref{table:complexity_static} of \cref{subsec:comparison_static} that the condition favoring model-based controllers is when
\begin{equation}\label{eq:bound_sensitivity_fast_accurate}
    \epsilon_{H,k} \triangleq \|\hat{H}_k - H_k\| \leq \frac{\bar{\epsilon}'}{(k+1)^\theta}, \quad \theta \geq \frac{1}{3}, \quad k \in \mathbb{N},
\end{equation}
where $\epsilon_{H,k}$ is the error of $\hat{H}_k$ relative to $H_k$, $\bar{\epsilon}' > 0$ is a pre-specified initial error bound, and $\theta$ is an exponent. For instance, \eqref{eq:bound_sensitivity_fast_accurate} holds when we use recursive least squares to learn sensitivities evolving by linear random processes\cite{picallo2021adaptive}. 
\revise{In this scenario, pursuing a gray-box pipeline may seem unnecessary. Nonetheless, for the sake of comparison and completeness, we present the associated coefficient design and analysis. The gray-box controller \eqref{eq:hybrid_controller} employs the following combination coefficient
\begin{equation}\label{eq:adaptive_comb_coeff_case_accurate}
    \alpha_k = 1 - \min\Big\{\frac{C}{p(k+1)^\frac{1}{3}}, 1\Big\}, \quad k \in \mathbb{N},
\end{equation}
where $C$ is a constant (e.g., $C \geq p$). As $k$ increases, $\alpha_k$ gradually changes from $0$ to $1$. The intuition is to begin with random exploration to facilitate excitation and online sensitivity learning, then progressively exploit refined sensitivity estimates for convergence.
}

However, different issues, e.g., noisy measurements, lack of covariance information, nonlinear dynamics, or resource constraints, can render sensitivity estimation slow or inaccurate, if not impossible. \revise{In these scenarios, our gray-box controller \eqref{eq:hybrid_controller} offers a principled approach of utilizing inaccurate sensitivities while retaining provable optimality guarantees.}
We distinguish two general cases related to $\hat{H}_k$ and present the corresponding strategies for tuning $\alpha_k$ in \eqref{eq:hybrid_controller_gradient}.



\textit{Case 1: Approximate Sensitivity with a Bounded Error} \par
In many applications, we construct approximate sensitivities based on prior knowledge or first-principles models\cite{ma2023reinforcement}, which are then fixed during online operation. This practice corresponds to the case where a bounded error exists between the approximate sensitivity $\hat{H}_k$ and the ground truth $H_k$, i.e.,
\begin{equation}\label{eq:bound_sensitivity_constant}
    \epsilon_{H,k} = \|\hat{H}_k - H_k\| \leq \bar{\epsilon}, \quad k \in \mathbb{N},
\end{equation}
where $\bar{\epsilon} > 0$.
In this case, we select a constant $C > 0$ and use the following vanishing combination coefficient
\begin{equation}\label{eq:adaptive_comb_coeff_case_constant}
    \alpha_k = \min\Big\{\frac{C p}{(k\!+\!1)^{\frac{1}{3}}}, 1 \Big\}, \quad k \in \mathbb{N}.
\end{equation}

\textit{Case 2: Asymptotically Accurate Sensitivity} \par
Online estimation techniques can be incorporated to generate increasingly accurate sensitivity estimates based on plant trajectories. However, we may not always learn sensitivities sufficiently fast due to measurement errors, lack of covariance data, etc. 
For instance, the estimation error may decrease as
\begin{equation}\label{eq:bound_sensitivity_decay}
    \epsilon_{H,k} = \|\hat{H}_k \!-\! H_k\| \leq \frac{\bar{\epsilon}'}{(k+1)^\theta}, \quad \theta \in \Big(0,\frac{1}{3}\Big), ~~ k \in \mathbb{N},
\end{equation}
where $\bar{\epsilon}' > 0$ is an error bound. In other words, the estimate $\hat{H}_k$ asymptotically converges to $H_k$, but the convergence rate is not sufficiently fast. In this case, we choose a constant $C>0$ and tune $\alpha_k$ according to
\begin{equation}\label{eq:adaptive_comb_coeff_case_decay}
    \alpha_k = \min\Big\{\frac{C p}{(k\!+\!1)^{\frac{1}{3}-\theta}},1 \Big\}, \quad k \in \mathbb{N}.
\end{equation}
We remark that the exponent $\theta$ in \eqref{eq:bound_sensitivity_decay} can be derived via modern non-asymptotic characterizations of learning systems. If an estimate of $\theta$ is inexact or elusive, we resort to the conservative scenario \eqref{eq:bound_sensitivity_constant} and select $\alpha_k$ according to \eqref{eq:adaptive_comb_coeff_case_constant}. 


Since $\alpha_k$ starts from $1$ and approaches $0$ as $k$ increases, the update rule \eqref{eq:hybrid_controller_gradient} integrating \eqref{eq:adaptive_comb_coeff_case_constant} or \eqref{eq:adaptive_comb_coeff_case_decay} initially favors the model-based inexact gradient $\tilde{\phi}_{k,1}$ for rapid response and later the model-free gradient estimate $\tilde{\phi}_{k,2}$ for solution accuracy. The combination coefficients allow smooth interpolation between model-based and model-free directions. This contributes to balanced performance despite variations in model quality. Further, sequential pipelines (first model-based, followed by model-free, or vice versa) are covered as special cases. For instance, given sensitivities with bounded errors as in \eqref{eq:bound_sensitivity_constant}, setting $\alpha_k=1$ for $k\leq \lceil \beta T^{\frac{1}{3}}\rceil$ and $\alpha_k=0$ thereafter indicates a pipeline that is initially model-based and then model-free. Overall, the adaptive combination \eqref{eq:hybrid_controller_gradient} endows our gray-box controller with flexibility in exploiting sensitivities of varying qualities.



\section{Performance Analysis and Comparison}\label{sec:analysis}
We establish performance certificates when our gray-box controller \eqref{eq:hybrid_controller} is interconnected in closed loop with the plant \eqref{eq:sys_map}. \revise{Our focus is to decouple the impacts of sensitivity errors (model-based) and variances of stochastic exploration (model-free) on the performance guarantee. We also compare with model-based and model-free feedback optimization controllers and discuss performance trade-offs.} 

\vspace{-2ex}
\subsection{Performance Certificates}\label{subsec:performance_ncvx}
We provide a recursive inequality for the second moment $\E_{v_{[k]}}[\|\tilde{\phi}_k\|^2]$ and highlight the roles of sensitivity errors and variances due to gradient estimation, where $\tilde{\phi}_k$ is the update direction, see \eqref{eq:hybrid_controller_GD_update}, and $\E_{v_{[k]}}[\cdot]$ denotes the expectation with respect to i.i.d.~samples $v_{[k]} \triangleq \{v_0,\ldots,v_{k}\}$. Recall that $\eta$, $\delta$, and $\alpha_k$ are the step size, the smoothing parameter, and the combination coefficient, respectively. Moreover, $p$ is the size of the input, $T$ is the number of iterations set beforehand, and $\epsilon_{H,k} = \|\hat{H}_k - H_k\|$ is the sensitivity error.

\begin{lemma}\label{lem:sec_mom_ubd}
    If \cref{assump:sys_map,assump:objective} hold, the update rule \eqref{eq:hybrid_controller} implies
    \begin{align}\label{eq:sec_mom_overall}
        &\E_{v_{[k]}}[\|\tilde{\phi}_k\|^2] \notag \\
        &\leq \frac{4(1\!-\!\alpha_k)^2}{\delta^2}p^2M^2\eta^2 \E_{v_{[k]}}[\|\tilde{\phi}_{k-1}\|^2] \!+\! 8\alpha_k^2 \E_{v_{[k]}}[\|\nabla\obj(w_k)\|^2] \notag \\
        &\quad + 48(1\!-\!\alpha_k)^2p^2 \E_{v_{[k]}}[\|\nabla\obj(w_{k-1})\|^2] + 4\alpha_k^2 M_\Phi^2 \epsilon_{H,k}^2 \notag \\
        &\quad  + 8\alpha_k^2L^2\delta^2 + 6(1\!-\!\alpha_k)^2L^2p^2\delta^2.
    \end{align}
\end{lemma}

In the upper bound \eqref{eq:sec_mom_overall}, the term $4\alpha_k^2 M_\Phi^2 \epsilon_{H,k}^2$ quantifies the effect of the sensitivity error, and \revise{the last two terms involving $\delta$ arise from the deviations and variances due to exploration noise in the gradient estimate}. The relative importance of these terms can be adjusted by $\alpha_k$. The inequality \eqref{eq:sec_mom_overall} reflects the progress of the controller and helps to establish our closed-loop performance certificate as follows.
\begin{theorem}\label{thm:optimality_general}
    \revise{
    Suppose that \cref{assump:sys_map,assump:objective,assump:objective_Lipschitz} hold. For any specified $T \in \mathbb{N}_{+}$, let $\delta = \sqrt{2}M/\big(112\max(C,1/C) L T^{\frac{1}{3}}\big)$ and $\eta$ satisfy
    \begin{equation*}
        0 < \eta \leq \frac{1}{224\max(C,\frac{1}{C}) L pT^\kappa}, \quad
        \kappa =
        \begin{cases}
            \frac{1}{3} & \text{if \eqref{eq:bound_sensitivity_fast_accurate} holds,} \\
            \frac{2}{3} & \text{if \eqref{eq:bound_sensitivity_constant} holds,} \\
            \frac{2}{3}-\theta & \text{if \eqref{eq:bound_sensitivity_decay} holds.}
        \end{cases}
    \end{equation*}
    The closed-loop interconnection of the plant \eqref{eq:sys_map} and the gray-box controller \eqref{eq:hybrid_controller} ensures
    \begin{align}\label{eq:convergence_measure}
        \frac{1}{T}&\sum_{k=0}^{T-1}\E_{v_{[T]}}[\|\nabla\obj(w_k)\|^2] \leq \frac{4(\obj(w_0) \!-\! \obj^*)}{\eta T} \!+\! \frac{4M_\Phi^2}{T} \sum_{k=0}^{T-1} \alpha_k \epsilon_{H,k}^2 \notag \\
        &\hspace{3em} + \frac{L^2 \delta^2}{2T} \sum_{k=0}^{T-1} \left(8\alpha_k + p^2(1\!-\!\alpha_k)\right) + \frac{C_{M,\|\tilde{\phi}_0\|^2}}{T},
    \end{align}
    where $C_{M,\|\tilde{\phi}_0\|^2} \!>\! 0$ is a constant that only depends on $M$ and $\|\tilde{\phi}_0\|^2$. 
    Further, if $\{\hat{H}_k\}_{k=0}^{T-1}$ satisfy \eqref{eq:bound_sensitivity_fast_accurate}, \eqref{eq:bound_sensitivity_constant}, or \eqref{eq:bound_sensitivity_decay}, $\{\alpha_k\}_{k=0}^{T-1}$ are designed by \eqref{eq:adaptive_comb_coeff_case_accurate}, \eqref{eq:adaptive_comb_coeff_case_constant}, \eqref{eq:adaptive_comb_coeff_case_decay}, respectively, and $\eta$ attains its upper bound, then the closed-loop interconnection guarantees
    }
    \begin{equation}\label{eq:convergence_measure_complexity}
        \frac{1}{T}\sum_{k=0}^{T-1}\E_{v_{[T]}}[\|\nabla\obj(w_k)\|^2] \!=\! 
        \begin{cases}
            \hspace{-3pt}\mathcal{O}\Big(\frac{p}{T^\frac{2}{3}}\Big) &\hspace{-4pt} \text{if \eqref{eq:bound_sensitivity_fast_accurate},\,\eqref{eq:adaptive_comb_coeff_case_accurate} hold,} \\
            \hspace{-3pt}\mathcal{O}\Big(\frac{p}{T^{\frac{1}{3}}}\Big) &\hspace{-4pt} \text{if \eqref{eq:bound_sensitivity_constant},\,\eqref{eq:adaptive_comb_coeff_case_constant} hold,} \\
            \hspace{-3pt}\mathcal{O}\Big(\frac{p}{T^{\frac{1}{3}+\theta}}\Big) &\hspace{-4pt} \text{if \eqref{eq:bound_sensitivity_decay},\,\eqref{eq:adaptive_comb_coeff_case_decay} hold.}
        \end{cases}
    \end{equation}
\end{theorem}

The performance measure in \cref{thm:optimality_general} is the average squared gradient norm of the reduced objective $\obj$ evaluated at the candidate solutions $\{w_k\}_{k=0}^{T-1}$, quantifying the local optimality of solutions. This measure is typical for nonconvex optimization (see \cite{nesterov2017random,tang2023zeroth}), given that achieving global optimality is NP-hard and intractable. 
Given any fixed a priori $T$, practitioners can choose constant $\eta$, $\delta$ and retain the above bounds up to time $T$. Afterward, the system performance does not degrade, although further improvement is no longer guaranteed. Alternatively, one may work with consecutive periods of a prescribed long length $T$. In each period, the controller \eqref{eq:hybrid_controller} is implemented in closed loop with the plant \eqref{eq:sys_map} with constant $\eta$, $\delta$ dependent on $T$, and the last solution of the previous period is used as a warm start. This scheme extends the presented non-asymptotic closed-loop certificates to scenarios without finite termination times.


The upper bound \eqref{eq:convergence_measure} elucidates how the model-based direction \eqref{eq:hybrid_controller_inexact_GD} and the model-free estimate \eqref{eq:hybrid_controller_grad_est} shape the overall closed-loop performance. Specifically, the term $\sum_{k=0}^{T\!-\!1} \alpha_k \epsilon_{H,k}^2 /T$ quantifies the cumulative error when the approximate sensitivity $\hat{H}_k$ is used instead of the true sensitivity $H_k$, whereas the second-to-last term is due to the deviations and variances arising from stochastic exploration in \eqref{eq:hybrid_controller_exploration}.
\revise{The combination coefficients $\{\alpha_k\}_{k=0}^{T-1}$ prescribed in \cref{thm:optimality_general} enable our gray-box controller \eqref{eq:hybrid_controller} to exploit approximate sensitivities of varying qualities, thus achieving balanced closed-loop performance in terms of sample efficiency, dimensional dependence, and provable optimality quantified by \eqref{eq:convergence_measure_complexity}.
}



When our gray-box controller \eqref{eq:hybrid_controller} is applied to stable nonlinear dynamical systems, we can establish performance guarantees by introducing converse Lyapunov functions and analyzing the coupled evolution of system dynamics and controller iterations, see, e.g., \cite{belgioioso2022online,he2022model}. Further, if the nonconvex problem~\eqref{eq:opt_original} involves an input constraint set, similar certificates can be derived by examining the second moment of the gradient mapping, see e.g.~\cite[Eq.~(2.4)]{ghadimi2016mini}.

\vspace{-1ex}
\subsection{Comparison with Model-based and Model-Free Controllers}\label{subsec:comparison_static}
\revise{Our gray-box controller \eqref{eq:hybrid_controller} features principled interpolation based on approximate sensitivities and stochastic exploration, thereby achieving balanced closed-loop performance. The unified formulation of design and analysis incorporates basic model-based and model-free pipelines as special cases. We establish the corresponding performance guarantees and demonstrate the trade-off between iteration and dimensional dependence.
}

We characterize the closed-loop performance incurred by the deterministic model-based controller \eqref{eq:model-based-FO} using approximate sensitivities $\{\hat{H}_k\}_{k=0}^{T-1}$ \revise{(or equivalently, \eqref{eq:hybrid_controller} with $\alpha_k=1, \forall k \in \mathbb{N}$ and $\delta=0$)}. 



\begin{theorem}\label{thm:model_based_optimality}
    Suppose that \cref{assump:sys_map,assump:objective,assump:objective_Lipschitz} hold. Consider the deterministic model-based controller whose update rule is \eqref{eq:model-based-FO} with $\nabla_u h(u_k,d)$ replaced by $\hat{H}_k$. 
    Let $0 < \eta \leq 1/(4L)$. The closed-loop interconnection of the plant \eqref{eq:sys_map} and this controller results in
    \begin{equation}\label{eq:convergence_measure_model_based}
        \frac{1}{T} \sum_{k=0}^{T-1}\|\nabla\obj(w_k)\|^2 \leq \revise{\frac{4\big(\obj(w_0) \!-\! \obj^*\big)}{\eta T}} + \frac{3M_\Phi^2}{T} \sum_{k=0}^{T-1} \epsilon_{H,k}^2.
    \end{equation}
    \revise{Moreover, when $\eta$ attains its upper bound, the order of the right-hand side of \eqref{eq:convergence_measure_model_based} depends on the quality of $\{\hat{H}_k\}_{k=0}^{T-1}$, as given by}
    \begin{equation}\label{eq:convergence_measure_complexity_MB}
        \frac{1}{T} \sum_{k=0}^{T-1}\|\nabla\obj(w_k)\|^2 \!=\!
        \begin{cases}
            \mathcal{O}\Big(\frac{1}{T^{\min(1,2\theta)}}\Big) \hspace{-1ex} & \text{if \eqref{eq:bound_sensitivity_fast_accurate} holds, } \theta \!\neq\! \frac{1}{2}, \\
            \mathcal{O}\Big(\frac{\ln T}{T}\Big) & \text{if \eqref{eq:bound_sensitivity_fast_accurate} holds, } \theta \!=\! \frac{1}{2}, \\
            \mathcal{O}\Big(\frac{1}{T} + \bar{\epsilon}^2\Big) & \text{if \eqref{eq:bound_sensitivity_constant} holds}, \\
            \bigO{\frac{1}{T^{2\theta}}} & \text{if \eqref{eq:bound_sensitivity_decay} holds}. \\
        \end{cases}
    \end{equation}
\end{theorem}


Further, we quantify the performance of the model-free controller \cite{he2022model}. This controller is a special case of \eqref{eq:hybrid_controller} with the combination coefficients $\alpha_k = 0, \forall k \in \mathbb{N}$, because only the gradient estimate \eqref{eq:hybrid_controller_grad_est} rather than $\hat{H}_k$ is used during iterations. 

\begin{theorem}\label{thm:optimality_general_mf}
    \revise{Let \cref{assump:sys_map,assump:objective,assump:objective_Lipschitz} hold, $\delta=\sqrt{2}M/(112LT^\frac{1}{3})$, and $0< \eta \leq 1/(224Lp^2T^\frac{1}{3})$. When $\eta$ attains its upper bound}, the closed-loop interconnection of the plant~\eqref{eq:sys_map} and the model-free controller (i.e., \eqref{eq:hybrid_controller} with $\alpha_k = 0, \forall k \in \mathbb{N}$) results in
    \begin{align*}
        &\frac{1}{T}\sum_{k=0}^{T-1}\E_{v_{[T]}}[\|\nabla\obj(w_k)\|^2] \notag \\
        &\leq \Big(\!896L(\obj(w_0) \!-\! \obj^*) \!+\! \frac{M^2}{112^2} \!\Big) \frac{p^2}{T^\frac{2}{3}} \!+\! \frac{C'_{M,\|\tilde{\phi}_0\|^2}}{T} \!=\! \bigO{\frac{p^2}{T^\frac{2}{3}}},
    \end{align*}
    where $C'_{M,\|\tilde{\phi}_0\|^2} > 0$ is a constant that only depends on $M$ and $\|\tilde{\phi}_0\|^2$ but not $p$ and $T$.
\end{theorem}


We summarize the performance guarantees of the model-based \eqref{eq:model-based-FO}, model-free (i.e., \eqref{eq:hybrid_controller} with $\alpha_k\!=\!0, \forall k$), and gray-box controllers for problem \eqref{eq:opt_original} in \cref{table:complexity_static}. The certificates of the model-based controller \eqref{eq:model-based-FO} are independent of the problem dimension $p$. 
If the sensitivity $\hat{H}_k$ is highly accurate, then the model-based controller achieves the most favorable performance among the three. However, if the accuracy condition \eqref{eq:bound_sensitivity_fast_accurate} is not met, cumulative errors $\sum_k \epsilon_{H,k}^2$ dominate in the bound \eqref{eq:convergence_measure_model_based}, degrading the overall convergence rate. The model-free controller attains the same $\mathcal{O}(p^2/T^{\frac{2}{3}})$ scaling across all the cases, since it does not benefit from approximate sensitivities. Nonetheless, there is a quadratic dependence on $p$, implying that the model-free controller may be inefficient for higher-dimensional problems.

In contrast, our gray-box controller \eqref{eq:hybrid_controller} ensures provable optimality despite sensitivity errors and preserves sample efficiency, achieving a balanced trade-off in its dependence on the problem dimension $p$ and the iteration count $T$. When $\hat{H}_k$ satisfies \eqref{eq:bound_sensitivity_constant} or \eqref{eq:bound_sensitivity_decay}, the gray-box controller achieves better scaling in $T$ than the model-based controller, but worse scaling than the model-free controller. Further, the gray-box controller incurs convergence bounds linearly dependent on $p$, improving the quadratic dependence of the model-free approach. 
Under \eqref{eq:bound_sensitivity_decay}, the rate of the model-free controller outperforms that of the gray-box controller only when $T > p^{\frac{3}{1-3\theta}}$ is sufficiently large, a requirement that becomes prohibitive when $\theta$ is close to $\frac{1}{3}$ or $p$ is large. In contrast, when \eqref{eq:bound_sensitivity_fast_accurate} is satisfied, the gray-box controller achieves more favorable scaling in both the iteration count $T$ and the problem dimension $p$ than the model-free approach. The balanced behavior of the gray-box controller arises from the principled interpolation offered by the combination coefficients.

\begin{table}[!tb]
\centering
\renewcommand \arraystretch{1}
\caption{Orders of convergence measures for problem \eqref{eq:opt_original}}
\label{table:complexity_static}
\begin{threeparttable}
    \begin{tabularx}{\columnwidth}{*{4}{Y}}
        \toprule
        \makecell{{\bfseries Controllers}} & \makecell{Highly \\ accurate \eqref{eq:bound_sensitivity_fast_accurate} \\ $\theta \geq \frac{1}{3}$} & \makecell{Asymptotically \\ accurate \eqref{eq:bound_sensitivity_decay} \\ $\theta \in (0, \frac{1}{3})$} & \makecell{Bounded \\ errors \eqref{eq:bound_sensitivity_constant}} \\
        \midrule
        \makecell{model-based \\ \eqref{eq:model-based-FO} with $\hat{H}_k$} & $\widetilde{\mathcal{O}}\Big(\frac{1}{T^{\min(1,2\theta)}}\Big)$\tnote{i} & $\bigO{\frac{1}{T^{2\theta}}}$ & $\bigO{\frac{1}{T}+\bar{\epsilon}^2}$ \\
        \midrule
        \makecell{model-free \eqref{eq:hybrid_controller} \\ with $\alpha_k=0$} & $\bigO{\frac{p^2}{T^{\frac{2}{3}}}}$ & $\bigO{\frac{p^2}{T^{\frac{2}{3}}}}$ & $\bigO{\frac{p^2}{T^{\frac{2}{3}}}}$ \\
        \midrule
        gray-box \eqref{eq:hybrid_controller} & $\mathcal{O}\Big(\frac{p}{T^\frac{2}{3}}\Big)$ & $\mathcal{O}\Big(\frac{p}{T^{\theta + \frac{1}{3}}}\Big)$ & $\mathcal{O}\Big(\frac{p}{T^{\frac{1}{3}}}\Big)$ \\
        \bottomrule
    \end{tabularx}
    \begin{tablenotes}[para]\footnotesize
        \item[i] We omit the logarithmic factor for $\theta=\frac{1}{2}$, see \eqref{eq:convergence_measure_complexity_MB}.
    \end{tablenotes}
\end{threeparttable}
\end{table}




\section{Constrained Time-Varying Feedback Optimization}\label{sec:tracking}
We explore regulating a plant given time-varying objective functions and exogenous disturbances subject to input constraints. In power systems applications, the objective shifts due to tracking time-varying set points for voltages or power flows, and the changing disturbance arises from volatile renewable generation\cite{simonetto2020time}. Moreover, control inputs are constrained because of physical actuation limits or coupled economic requirements. \revise{This setting requires controller design and analysis that address the variation of optimal solutions. While there is abundant work in online optimization and control (see \cref{subsec:TV_performance}), our agenda is to explicitly decouple how sensitivity errors (model-based) and variances of stochastic exploration (model-free) impact online performance relative to the optimal benchmark.}

\subsection{Problem Formulation}\label{subsec:formulation-TV}
The specifications of changing objectives, variable disturbances, and input constraints lead to a constrained time-varying problem:
\begin{equation}\label{eq:opt_TV}
\begin{split}
    \min_{u\in \mathbb{R}^p, y\in \mathbb{R}^q} ~ & \Phi_k(u,y) \\ 
    \text{s.t.} ~ & y = h(u,d_k), \\
                ~& u \in \mathcal{U},
\end{split}
\end{equation}
where $\Phi_k:\mathbb{R}^p\times \mathbb{R}^q \to \mathbb{R}$ and $d_k \in \mathbb{R}^r$ are the objective and the unknown disturbance at time $k \in \mathbb{N}$, respectively. 
Moreover, $\mathcal{U} \subset \mathbb{R}^p$ is the constraint set for $u$. Problem~\eqref{eq:opt_TV} involves optimizing the input-output performance of the following time-varying steady-state map
\begin{equation}\label{eq:sys_map_TV}
    y = h(u,d_k).
\end{equation}
\revise{In the above discrete-time formulation, at each time interval $[k, k+1]$, we iteratively adjust the input $u$ after the output $y$ encoding $d_k$ is measured and the current objective $\Phi_k$ is revealed. While $\{d_k\}_{k \in \mathbb{N}}$ vary across time steps, they are processed as constants within each interval. This setup, for example, arises in continuous-time plants subject to changing disturbances and controlled in the sampled-data framework with zero-order holds \cite{belgioioso2022online}.} 
Let $\obj_k(u) \triangleq \Phi_k(u,h(u,d_k))$ be the reduced objective at time $k$. Our assumptions are as follows. 

\begin{assumption}\label{assump:constraint_set}
    The set $\mathcal{U}$ is a compact and convex set with diameter $D>0$, i.e., $\forall u_1,u_2 \in \mathcal{U}, \|u_1-u_2\|\leq D$. 
\end{assumption}

\begin{assumption}\label{assump:objective_TV_cvx}
    The function $\obj_k(u)$ is convex, $L_k$-smooth, and $M_k$-Lipschitz with respect to $u$. The function $\Phi_k(u,y)$ is $M_{\Phi,k}$-Lipschitz with respect to $y$. Moreover, $\{L_k\}$, $\{M_k\}$, and $\{M_{\Phi,k}\}$ are bounded.
\end{assumption}

Let $\mathcal{U}_{\tau} \triangleq \{u+\tau v \mid u\in \mathcal{U},v\in \mathbb{B}_p\}$ denote a set inflated from $\mathcal{U}$ by a limited range $\tau \mathbb{B}_p$, where $\tau > 0$ is an expansion coefficient, and $\mathbb{B}_p$ is the closed unit ball in $\mathbb{R}^p$. The following assumption specifies the boundedness of $\obj_k$. 

\begin{assumption}\label{assump:objective_uniform_bound}
    The function $\obj_k(u)$ is uniformly bounded, i.e., $\exists G \geq 0, \exists \tau>0, \forall u \in \mathcal{U}_{\tau}, \forall k \in \mathbb{N}, |\obj_k(u)| \leq G$.
\end{assumption}


\cref{assump:constraint_set,assump:objective_TV_cvx,assump:objective_uniform_bound} are typical in the literature, e.g., \cite{zhao2021bandit,hazan2022introduction}. \cref{assump:constraint_set} implies that the norm of any point in $\mathcal{U}$ is bounded, i.e., $\exists \bar{D} \geq 0, \forall u \in \mathcal{U}, \|u\| \leq \bar{D}$. \cref{assump:objective_uniform_bound} is related to \cref{assump:constraint_set,assump:objective_TV_cvx}, because a continuous function defined on a compact set is bounded. 
We introduce the inflated set $\mathcal{U}_\tau$ to handle the scenario where the input, perturbed by exploration noise, may lie outside of the constraint set $\mathcal{U}$. In practice, if the objective $\obj_k$ is bounded on $\mathcal{U}$ and Lipschitz continuous, then it is also bounded on $\mathcal{U}_\tau$. 
Here, we consider time-varying convex objectives to offer clear and accessible characterizations of the adaptation performance through the corresponding guarantees. 
The analysis can be extended to tackle nonconvex objectives, although more involved analysis is required under specific performance metrics, e.g., the cumulative gradient norms of window-smoothed objectives \cite{hazan2017efficient}.

The unknown map $h$ and the changing disturbances $\{d_k\}$ prevent us from solving \eqref{eq:opt_TV} directly via numerical solvers. In contrast, 
we aim for an online closed-loop strategy, featuring a feedback optimization controller that exploits output measurements to optimize the dynamic behavior of \eqref{eq:sys_map_TV}. Let $u_k^* \in \mathbb{R}^p$ be an optimal point of problem~\eqref{eq:opt_TV} at time $k$. The goal is to generate control inputs that are competitive with the sequence of optimal solutions $\{u_k^*\}_{k\in \mathbb{N}}$.


\subsection{Design of the Running Gray-Box Controller}
To handle the constrained and time-varying problem~\eqref{eq:opt_TV}, we adjust our gray-box controller \eqref{eq:hybrid_controller} by leveraging projection and the most recent output measurement. The update rules are
\begin{subequations}\label{eq:hybrid_controller_TV}
\begin{align}
    w_{k+1} &= \proj(w_k - \eta \hat{\phi}_k), \label{eq:hybrid_controller_TV_GD_update} \\
    \hat{\phi}_k &= \alpha_k \hat{\phi}_{k,1} + (1-\alpha_k) \hat{\phi}_{k,2}, \label{eq:hybrid_gradient_TV} \\
    \hat{\phi}_{k,1} &= \nabla_{u} \Phi_k(u_k,y_k) + \hat{H}_k \nabla_y \Phi_k(u_k,y_k), \label{eq:hybrid_gradient_TV_MB} \\
    \hat{\phi}_{k,2} &= \frac{pv_k}{\delta} \big(\Phi_k(u_k,y_k) - \Phi_{k-1}(u_{k-1},y_{k-1})\big), \label{eq:hybrid_gradient_TV_MF} \\
    u_{k+1} &= w_{k+1} + \delta v_{k+1},  \label{eq:hybrid_TV_noise}
\end{align}
\end{subequations}
where $\proj(\cdot)$ denotes the projection to the constraint set $\mathcal{U}$ \revise{(defined as $\proj(u) \triangleq \argmin_{u' \in \mathcal{U}} \|u-u'\|$)}, and the involved variables are the same as \eqref{eq:hybrid_controller}, see also the discussions below \eqref{eq:hybrid_controller}. 
Analogously, in the iterative update, our running gray-box controller \eqref{eq:hybrid_controller_TV} merges the model-based inexact gradient \eqref{eq:hybrid_gradient_TV_MB} and the model-free gradient estimate \eqref{eq:hybrid_gradient_TV_MF}. Subsequently, projection to $\mathcal{U}$ is performed to ensure constraint satisfaction. While the descent-based updates of \eqref{eq:hybrid_controller} and \eqref{eq:hybrid_controller_TV} resemble each other, they entail different performance guarantees when applied to the corresponding problems.



\begin{remark}
    While the candidate solution $w_k$ lies in $\mathcal{U}$, the input $u_k$ in the transient stage may violate the constraint. If we need strict constraint satisfaction, we can project in \eqref{eq:hybrid_controller_TV_GD_update} onto a deflated set $(1\!-\!\kappa) \mathcal{U}$ as \cite{zhao2021bandit}, where $\kappa \in (0,1)$.
\end{remark}

Similar to \cref{subsec:adaptive_weight}, for problem~\eqref{eq:opt_TV}, model-based controllers using $\{\hat{H}_k\}_{k\in \mathbb{N}}$ (i.e., \eqref{eq:hybrid_controller_TV} with $\alpha_k \!=\! 1,\forall k\in \mathbb{N}$ and $\delta \!=\! 0$) are favorable \revise{if $\hat{H}_k$ is a highly accurate estimate of the sensitivity matrix $H_k' \triangleq \nabla_u h(u_k,d_k)$, namely,}
\begin{equation}\label{eq:bound_sens_fast_accurate_TV}
    \epsilon_{H,k} \triangleq \|\hat{H}_k - H_k'\| \leq \frac{\bar{\epsilon}'}{(k+1)^\theta}, \quad \theta \geq \frac{1}{3}, \quad k \in \mathbb{N},
\end{equation}
where $\epsilon_{H,k}$ is the error of the approximate sensitivity, and $\bar{\epsilon}' > 0$. \revise{For the sake of completeness, we use the following rule to tune $\alpha_k$:
\begin{equation}\label{eq:adaptive_coeff_TV_case_accurate}
    \alpha_k = 1 - \min \Big\{\frac{C}{p^{\frac{1}{6}}\sqrt{k+1}}, 1\Big\}.
\end{equation}
The intuition is to begin with stochastic exploration for excitation and learning, and gradually shift to model-based updates, since approximate sensitivities are generally reliable.
}

Nonetheless, \eqref{eq:bound_sens_fast_accurate_TV} may not always hold due to various issues, e.g., noisy measurements or nonlinear dynamics. \revise{We present the selection of $\{\alpha_k\}_{k \in \mathbb{N}}$ based on the quality of $\hat{H}_k$:
\begin{subnumcases}{
\alpha_k \!=\! \label{eq:adaptive_coeff_TV}}
    {\textstyle \min\Big\{\frac{C p^\frac{1}{6}}{(k+1)^{\frac{1}{6}}}, 1 \Big\}}, & if \eqref{eq:bound_sensitivity_constant} holds, \label{eq:bound_sensitivity_constant_TV} \\
    {\textstyle \min\Big\{\frac{C p^\frac{1}{6}}{(k+1)^{\max(\frac{1}{6}-\frac{\theta}{2}, \frac{1}{12})}}, 1\Big\}}, & if \eqref{eq:bound_sensitivity_decay} holds. \label{eq:bound_sens_decay_TV}
\end{subnumcases}
Applications where $\theta$ in \eqref{eq:bound_sens_decay_TV} is hard to establish can be treated through the more conservative case \eqref{eq:bound_sensitivity_constant_TV}.
}

\subsection{Performance Certificates}\label{subsec:TV_performance}

\revise{
We establish performance certificates for the closed-loop interconnection between the plant \eqref{eq:sys_map_TV} and our running gray-box controller \eqref{eq:hybrid_controller_TV}. Two primary non-asymptotic performance measures exist for the time-varying problem \eqref{eq:opt_TV}. The first is the finite-time tracking error \cite{simonetto2020time,ospina2022feedback}, which offers last-iterate guarantees in the spirit of input-to-state stability (see, e.g., \cite{belgioioso2022online,bianchin2021time,cothren2022online}). 
The second measure, which we address here, is the \emph{dynamic regret} capturing the cumulative performance gap relative to a sequence of benchmark solutions over a finite horizon \cite{zhao2021bandit,hazan2022introduction}.

Different from the above works, we explicitly characterize the influences of the sensitivity error $\epsilon_{H,k} \!=\! \|\hat{H}_k \!-\! H_k'\|$ and variances due to stochastic exploration on the performance measure. This characterization facilitates tuning the combination coefficients to achieve balanced closed-loop behaviors. We introduce a lemma that bounds the inner product and the second moment related to $\hat{\phi}_k$.

\begin{lemma}\label{lem:hybrid_prod_mom_bd}
    Let \cref{assump:sys_map,assump:constraint_set,assump:objective_TV_cvx,assump:objective_uniform_bound} hold. The rule \eqref{eq:hybrid_controller_TV} ensures
    \begin{subequations}
        \begin{align}
            &\E_{v_{[k]}}\big[\hat{\phi}_k^\top(w_k\!-\!u_k^*)\big] \geq \E_{v_{[k]}}[\obj_k(w_k)] \!-\! \obj_k(u_k^*) \notag \\
            &\hspace{8em} \!-\! \alpha_k D(M_{\Phi,k} \epsilon_{H,k} \!\!+\!\! L_k \delta) \!-\! (1\!\!-\!\!\alpha_k) L_k \delta^2, \label{eq:hybrid_prod_lbd} \\
            &\E_{v_{[k]}}\big[\|\hat{\phi}_k\|^2\big] \!\leq\! 4\alpha_k^2(M_k^2 \!\!+\!\! M_{\Phi,k}^2 \epsilon_{H,k}^2) \!+\! 8(1\!\!-\!\!\alpha_k)^2 \frac{p^2G^2}{\delta^2}. \label{eq:hybrid_mom_ubd}
        \end{align}
    \end{subequations}
\end{lemma}

In the bounds in \cref{lem:hybrid_prod_mom_bd}, the terms containing $\alpha_k$ and $1-\alpha_k$ quantify the contributions of the sensitivity error in the model-based direction \eqref{eq:hybrid_gradient_TV_MB} and the stochastic exploration required by the model-free estimate \eqref{eq:hybrid_gradient_TV_MF}, respectively. These joint influences allow systematic adjustment via $\alpha_k$ prescribed in \eqref{eq:adaptive_coeff_TV}.}


We characterize the online closed-loop performance through the dynamic regret $\reg \triangleq \sum_{k=1}^{T} \!\big(\E_{v_{[T]}}[\obj_k(w_k)] \!-\! \obj_k(u_k^*)\big)$, i.e., the cumulative difference between the expected objective values evaluated at the candidate solutions $\{w_k\}_{k=1}^{T}$ and at the optimal points $\{u_k^*\}_{k=1}^{T}$. To capture the variation of \eqref{eq:opt_TV}, we introduce the path length $C_T \triangleq \sum_{k=1}^{T} \|u_k^* \!-\! u_{k-1}^*\|$, which accumulates the shifts between two consecutive optimal points \cite{zhao2021bandit}. \revise{The resulting dynamic regret scales with the horizon $T$, the path length $C_T$, the dimension $p$ of the input, and the exponent $\theta$ in \eqref{eq:bound_sensitivity_fast_accurate} or \eqref{eq:bound_sensitivity_decay}.}



\begin{theorem}\label{thm:track_dynamic_regret}
    Suppose that \cref{assump:sys_map,assump:constraint_set,assump:objective_TV_cvx,assump:objective_uniform_bound} are satisfied. For any specified $T \in \mathbb{N}_{+}$, let $\eta = 1/\big(\sqrt{p}T^{\kappa}\big)$ and $\delta = \min\big(p^\frac{1}{3}/T^{\kappa'}, \tau\big)$, where $\kappa = \kappa' =\frac{1}{2}$ if \eqref{eq:bound_sensitivity_fast_accurate} holds; $\kappa=\frac{5}{6}, \kappa' = \frac{1}{6}$ if \eqref{eq:bound_sensitivity_constant} holds; $\kappa=\max(\frac{5}{6}-\frac{\theta}{2}, \frac{3}{4}), \kappa' = \frac{1}{6}$ if \eqref{eq:bound_sensitivity_decay} holds. 
    Given the conditions of $\{\hat{H}_k\}$ and the designs of $\{\alpha_k\}$ in \eqref{eq:adaptive_coeff_TV_case_accurate} and \eqref{eq:adaptive_coeff_TV}, the closed-loop interconnection of the plant \eqref{eq:sys_map_TV} and the gray-box controller \eqref{eq:hybrid_controller_TV} incurs the following dynamic regret bounds
    \begin{equation}\label{eq:dynamic_regret_complexity}
        \hspace{-3.5pt} \reg \!=\! 
        \begin{cases}
            \hspace{-3.5pt} \bigO{\!\sqrt{pT}(C_T\!+\!1) \!+\! T^{\max(\frac{1}{2},1\!-\!\theta)}\!} &\hspace{-6pt} \text{if \eqref{eq:bound_sensitivity_fast_accurate},\,\eqref{eq:adaptive_coeff_TV_case_accurate} hold,} \\
            \hspace{-3.5pt} \bigO{\sqrt{p} T^{\frac{5}{6}}(C_T \!+\! 1)} &\hspace{-6pt} \text{if \eqref{eq:bound_sensitivity_constant},\,\eqref{eq:bound_sensitivity_constant_TV} hold,} \\
            \hspace{-3pt} \bigO{\sqrt{p} T^{\max(\frac{5}{6}-\frac{\theta}{2},\frac{3}{4})}(C_T \!+\! 1)} &\hspace{-6pt} \text{if \eqref{eq:bound_sensitivity_decay},\,\eqref{eq:bound_sens_decay_TV} hold.}
        \end{cases}
    \end{equation}
\end{theorem}

\revise{The orders of the dynamic regret in \eqref{eq:dynamic_regret_complexity} feature a sublinear dependence on the problem dimension $p$. The scaling in the time horizon $T$ varies with the quality of the sensitivity. If $C_T$ is known in advance (e.g., via quantitative sensitivity bounds), then the dependence of $\reg$ on $C_T$ can be reduced from linear to square-root by choosing $\eta \sim \sqrt{C_T}/(\sqrt{p} T^\kappa)$. 
In practice, practitioners can 
leverage the scheme of consecutive periods discussed after \cref{thm:optimality_general} in \cref{subsec:performance_ncvx} to extend the above non-asymptotic certificates to scenarios without finite termination times.
Furthermore, in the time-varying setting where \eqref{eq:hybrid_controller_TV} is interconnected with a stable dynamical system, a promising means of characterizing convergence measures as \eqref{eq:dynamic_regret_complexity} is to leverage stability certificates (e.g., converse Lyapunov functions, input-to-state stability, or contraction metrics) and time-scale separation arguments \cite{belgioioso2022online}.
}

\revise{
The design and analysis of our running gray-box controller \eqref{eq:hybrid_controller_TV} encompass model-based and model-free feedback optimization controllers as special cases. The following theorem establishes the performance certificates of these controllers in the time-varying setting.
\begin{theorem}\label{thm:track_dynamic_regret_mf}
    Let \cref{assump:sys_map,assump:constraint_set,assump:objective_TV_cvx,assump:objective_uniform_bound} hold. The closed-loop interconnection of the plant \eqref{eq:sys_map_TV} and the model-based controller, i.e., \eqref{eq:hybrid_controller_TV} with $\alpha_k=1,\forall k\in \mathbb{N}$, $\eta=1/\sqrt{T}$, and $\delta=0$, ensures
    \begin{align}\label{eq:dynamic_regret_mb}
        \reg =
        \begin{cases}
            \bigO{\sqrt{T}C_T\!+\!T^{\max(\frac{1}{2},1-\theta)}} &\text{if \eqref{eq:bound_sensitivity_fast_accurate} holds,} \\
            \bigO{\sqrt{T}(C_T+1)+T} &\text{if \eqref{eq:bound_sensitivity_constant} holds,} \\
            \bigO{\sqrt{T}(C_T+1)+T^{1-\theta}} &\text{if \eqref{eq:bound_sensitivity_decay} holds}.
        \end{cases}
    \end{align}
    The interconnection of the plant \eqref{eq:sys_map_TV} and the model-free controller, i.e., \eqref{eq:hybrid_controller_TV} with $\alpha_k=0,\forall k\in \mathbb{N}$, $\eta = 1/(p^\frac{2}{3}T^\frac{2}{3})$, and $\delta = \min(p^\frac{1}{3}/T^\frac{1}{6},\tau)$, incurs $\reg = \bigO{p^\frac{2}{3}T^\frac{2}{3}(C_T \!+\! 1)}$.
\end{theorem}
}

\revise{
\ifarxiv
\cref{table:complexity_TV} in \cref{appendix:track_dynamic_regret_mf} 
\else
Table~II in \cite[Appendix~H]{he2024gray} 
\fi
summarizes the performance certificates of different controllers for problem~\eqref{eq:opt_TV}. The dynamic regret of the model-based controller (i.e., \eqref{eq:hybrid_controller_TV} with $\alpha_k=1$) is independent of the problem dimension. Nonetheless, when sensitivity matrices are of low quality (e.g., when \eqref{eq:bound_sensitivity_constant} holds), the cumulative sensitivity errors dominate the regret bound, causing a degradation with respect to the horizon $T$. In contrast, the model-free controller (i.e., \eqref{eq:hybrid_controller_TV} with $\alpha_k=0$) exhibits the same dynamic regret (with an $\mathcal{O}\big(p^\frac{2}{3}\big)$ scaling) in all cases, because its update does not involve sensitivity matrices.

Further comparisons are in order. In terms of iteration dependence, the dynamic regret incurred by the gray-box controller scales more favorably than that of the model-based controller when \eqref{eq:bound_sensitivity_constant} or \eqref{eq:bound_sensitivity_decay} with $\theta \in (0,\frac{1}{4})$ holds. Although this scaling is slightly worse than that of the model-free controller under these same conditions, the gray-box controller improves the iteration dependence when the sensitivity is highly accurate as per \eqref{eq:bound_sensitivity_fast_accurate}. Regarding dimensional dependence, the gray-box controller lies in the middle: it achieves an $\bigO{\sqrt{p}}$ scaling better (worse) than the model-free (model-based, respectively) controller. This balanced trade-off between iteration and dimensional dependence resembles the behavior in \cref{subsec:comparison_static} and arises from the principled interpolation \eqref{eq:hybrid_gradient_TV}.
}

\section{Numerical Evaluations}\label{sec:experiment}
We numerically evaluate the performance of our gray-box controllers. We specify all the parameters related to the problems and controllers in 
\ifarxiv
\cref{appendix:experiment_details}. Our code is available \cite{he2024grayboxcode}.
\else
Appendix~I of the online report \cite{he2024gray}. Our code is available at \url{https://github.com/zyhe/Gray-box-nonlinear-feedback-optimization.git}.
\fi

We establish theoretical results for algebraic maps. Nonetheless, as discussed in \cref{subsec:formulation}, the online feedback nature of our controllers enables performance optimization of stable dynamical systems. We pursue steady-state optimization of a system given by
\begin{equation}\label{eq:sys_simulation}
    \begin{aligned}
        x_{k+1} &= Ax_k + B_1 u_k + B_2 \left(\sin(u_k) + u_k^2\right) + Ed_x, \\
        y_k &= Cx_k + Dd_y,
    \end{aligned}
\end{equation}
where $x\in \mathbb{R}^{30}$, $u\in \mathbb{R}^{15}$, and $y\in \mathbb{R}^{10}$ denote the state, input, and output, respectively, and $d_x, d_y\in \mathbb{R}^{10}$ are disturbances.

First, we focus on the unconstrained nonconvex problem
\begin{align}\label{eq:opt_problem_simulation}
    &\min_{u \in \mathbb{R}^{15}, y \in \mathbb{R}^{10}} ~ \Phi(u,y)= - \lambda \|u\|^3 + u^{\top}M_1u + m_2^{\top}u + \|y\|^2 \notag \\
    &\hspace{2em} \text{s.t.} ~ y = C(I\!-\!A)^{-1} \left(B_1 u \!+\! B_2 (\sin(u) + u^2) \!+\! Ed_x\right) \!+\! Dd_y \notag \\
    &\hspace{2em} \phantom{\text{s.t.} ~ y} \triangleq h'(u,d_x,d_y).
\end{align}
Due to $-\lambda \|u\|^3$ and the nonlinear part $\sin(u) \!+\! u^2$, the reduced objective $\obj(u) \triangleq \Phi(u,h'(u,d_x,d_y))$ is nonconvex. Moreover, the equality constraint in \eqref{eq:opt_problem_simulation} corresponds to the steady-state map of the plant~\eqref{eq:sys_simulation}. The input-output sensitivity $H_k$ of \eqref{eq:sys_simulation} at $u_k$ is $H_k = \left[C(I \!-\! A)^{-1}(B_1 + B_2\diag(\cos(u_k) \!+\! 2 u_k))\right]^\top$,
where $\diag(\cos(u_k) \!+\! 2 u_k)$ is a diagonal matrix. 

We compare the closed-loop interconnection of \eqref{eq:sys_simulation} with various controllers: model-based feedback optimization controllers with $H_k$, inexact $\hat{H}$, and sensitivity learning \cite{picallo2021adaptive}, the discrete-time stochastic extremum seeking controller \cite{liu2016stochastic}, the model-free controller \cite{he2022model}, and our gray-box controller \eqref{eq:hybrid_controller} using $\hat{H}$. Since the last three controllers involve stochasticity, we conduct $30$ independent experiments.

\begin{figure}[t]
\centering
  \vspace{-2ex}
  \subfloat[Performance on problem~\eqref{eq:opt_problem_simulation}]{\includegraphics[width=0.5\columnwidth]{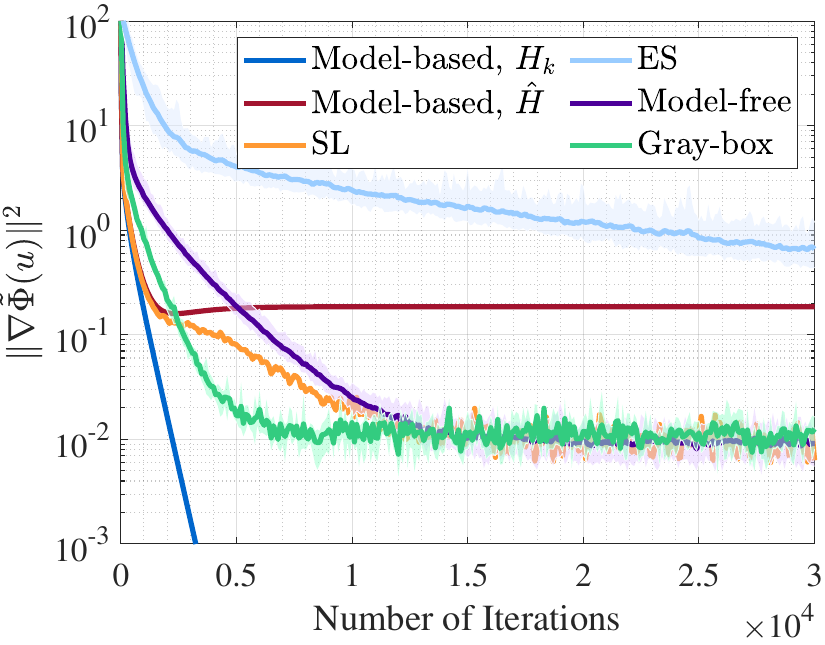}\label{fig:comparison}} \hfil
  \subfloat[Performance on problem~\eqref{eq:opt_problem_simulation_TV}]{\includegraphics[width=0.5\columnwidth]{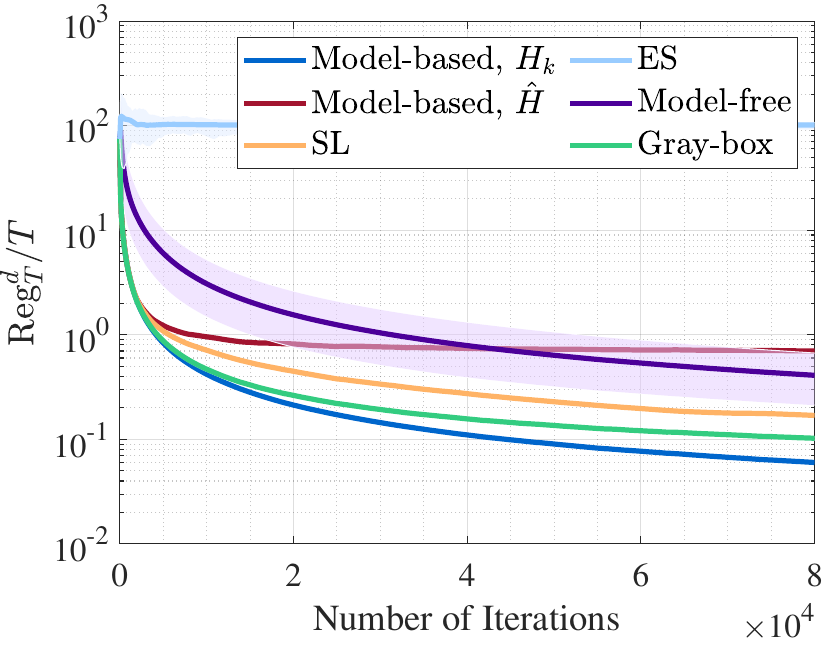}\label{fig:comparison_dynamic_regret}}
  \caption{\revise{Comparison of different controllers for problems \eqref{eq:opt_problem_simulation} and \eqref{eq:opt_problem_simulation_TV}. ``SL'' and ``ES'' denote ``sensitivity learning'' and ``Extremum Seeking'', respectively.}}
  \label{fig:convergence_curves}
\end{figure}

\cref{fig:comparison} illustrates the performance of the closed-loop interconnection. The convergence measure is the squared gradient norm of the reduced objective $\obj$ of \eqref{eq:opt_problem_simulation}. The shaded regions represent the variability across independent experiments, and the solid curves indicate the average convergence measures. When the true sensitivity $H_k$ is available, the model-based controller enjoys both fast convergence and high accuracy. However, the use of an inexact sensitivity $\hat{H}$ causes a severe bias and closed-loop sub-optimality. Sensitivity learning addresses this issue via an additional recursive estimation module. This module incurs additional costs of storage and computation. The extremum-seeking controller suffers from slow convergence. The model-free controller yields rather accurate solutions, though there is an increase in the required number of iterations. In contrast, our gray-box controller strikes a balance between the convergence rate and solution accuracy. Furthermore, the gray-box controller is easy to implement, in that it merely incorporates $\hat{H}$ via adaptive convex combination without requiring more accurate sensitivity estimates.


Next, we consider time-varying performance optimization of the plant \eqref{eq:sys_simulation} with input constraints:
\begin{equation}\label{eq:opt_problem_simulation_TV}
\begin{split}
    \min_{u \in \mathbb{R}^{15}, y \in \mathbb{R}^{10}} &~ \Phi_k(u,y)= u^{\top} M_{1,k} u + m_{2,k}^{\top} u + \|y\|^2 \\
    \text{s.t.} &~ y = h'(u,d_{x,k},d_{y,k}), \\
    &~ \underline{u} \leq u \leq \bar{u},
\end{split}
\end{equation}
where $h'(u,d_{x,k},d_{y,k})$ is the steady-state map of the plant~\eqref{eq:sys_simulation}, and $\underline{u} \in \mathbb{R}^{15}$ and $\bar{u} \in \mathbb{R}^{15}$ denote the lower bound and the upper bound on $u$, respectively. To align with the theoretical results in \cref{sec:tracking}, we consider the case where $B_2 = 0$ (i.e., the plant~\eqref{eq:sys_simulation} is linear), so that the overall objective of \eqref{eq:opt_problem_simulation_TV} is a convex function of $u$.


We augment the above controllers with projection onto the constraint set (similar to \eqref{eq:hybrid_controller_TV_GD_update}) and implement them in closed loop with the plant~\eqref{eq:sys_simulation}. Analogous to the above experiment, we perform $30$ independent runs.
\cref{fig:comparison_dynamic_regret} illustrates the evolutions of the time-averaged dynamic regret values (i.e., $\reg/T$) incurred by these closed-loop interconnections. We observe similar patterns as \cref{fig:comparison}. The direct use of the approximate sensitivity $\hat{H}$ diminishes solution accuracy. Nonetheless, by suitably incorporating this information, the gray-box controller achieves better performance compared to the model-free controller and the controller with sensitivity learning. Furthermore, for the considered iteration range, the gray-box controller is the closest to the benchmark with the exact sensitivity.



\section{Conclusion}\label{sec:conclusion}
We proposed gray-box feedback optimization controllers to optimize the steady-state performance of a nonlinear system in closed loop. \revise{These controllers merge approximate input-output sensitivities of the system into model-free updates via a tunable convex combination. Based on a principled interpolation between model-based and model-free approaches, we provided unified performance characterizations encompassing different approaches for unconstrained static and constrained time-varying problems. 
The proposed controllers achieve balanced closed-loop performance in sample efficiency (in terms of iteration and dimensional dependence) and provable accuracy, even in the presence of sensitivity errors.}
Future directions include leveraging other forms of model information, tackling output constraints via dualization, analyzing the interplay between model-free control and online identification, and \revise{exploring gray-box pipelines in the broader landscape of adaptive control}.

    \ifarxiv
        \balance
        \bibliographystyle{IEEEtran}
        \bibliography{article}
        
        \balance
\appendix
\restoreIEEEdisplayskips
\crefalias{section}{appendix}
\crefalias{subsection}{appendix}

\subsection{Auxiliary Lemmas}\label{appendix:cross_expt}
    \revise{Let the error between the true gradient $\nabla \obj(u_k)$ and the inexact gradient $\tilde{\phi}_{k,1}$ \eqref{eq:hybrid_controller_inexact_GD} using $\hat{H}_k$ be given by}
    \begin{align}\label{eq:grad_est_error}
        \epsilon_k \!\triangleq\! \nabla\obj(u_k) \!-\! \tilde{\phi}_{k,1} \!=\! (H_k \!-\! \hat{H}_k) \nabla_y \Phi(u_k,y_k), \quad k \!\in\! \mathbb{N}.
    \end{align}
    The following lemma provides an upper bound on the cross term $\E_{v_{[k]}}[-\nabla \obj(w_k)^{\top}\tilde{\phi}_k]$, where $\nabla \obj(w_k)$ is the gradient of the objective function at the candidate solution $w_k$, and $-\tilde{\phi}_k$ is the update direction of the controller, see \eqref{eq:hybrid_controller_gradient}. This bound will be useful for establishing closed-loop performance certificates.
    \begin{lemma}\label{lem:cross_expt_ubd}
        If \cref{assump:sys_map,assump:objective} hold, with \eqref{eq:hybrid_controller}, we have
        \begin{align}\label{eq:cross_expt_overall}
            \E_{v_{[k]}}&[-\nabla \obj(w_k)^{\top}\tilde{\phi}_k] \notag \\
                &\leq -\frac{1}{2}\E_{v_{[k]}}[\|\nabla\obj(w_k)\|^2] + \alpha_k\E_{v_{[k]}}[\|\epsilon_k\|^2] \notag \\
                &\quad + \alpha_kL^2\delta^2 + \frac{1-\alpha_k}{8}L^2p^2\delta^2.
        \end{align}
    \end{lemma}

    \begin{proof}
    The cross term $\E_{v_{[k]}}[-\nabla \obj(w_k)^{\top}\tilde{\phi}_k]$ satisfies
    \begin{align}\label{eq:expectation_mid_term}
        \E_{v_{[k]}}&[-\nabla \obj(w_k)^{\top}\tilde{\phi}_k] \notag \\
            &\stackrel{\text{(s.1)}}{=} \E_{v_{[k-1]}}\big[\E_{v_k}[-\nabla \obj(w_k)^{\top}\tilde{\phi}_k|v_{[k-1]}]\big] \notag \\
            &\stackrel{\text{(s.2)}}{=} \E_{v_{[k-1]}}\big[\nabla \obj(w_k)^{\top}\E_{v_k}[-\tilde{\phi}_k|v_{[k-1]}]\big],
    \end{align}
    where (s.1) utilizes the tower rule, and (s.2) holds because $\obj(w_k)$ is measurable with respect to $v_{[k-1]}$.

    Let $\obj_{\delta}: \mathbb{R}^p \to \mathbb{R}$ be the smooth approximation of the objective $\obj$, see also the definition \eqref{eq:smooth_approx}.
    For the term $\E_{v_k}[-\tilde{\phi}_k|v_{[k-1]}]$, we have
    \begin{align}\label{eq:cond_expt_grad_est}
        \E_{v_k}&[-\tilde{\phi}_k|v_{[k-1]}] \notag \\
            &= -\alpha_k \E_{v_k}[\tilde{\phi}_{k,1}|v_{[k-1]}] - (1-\alpha_k) \E_{v_k}[\tilde{\phi}_{k,2}|v_{[k-1]}] \notag \\
            &\stackrel{\text{(s.1)}}{=} -\alpha_k\nabla\obj(w_k) + \alpha_k\E_{v_k}[\nabla\obj(w_k) - \tilde{\phi}_{k,1} | v_{[k-1]}] \notag \\
            &\quad - (1-\alpha_k)\nabla\obj_{\delta}(w_k) \notag \\
            &\stackrel{\text{(s.2)}}{=} -\nabla\obj(w_k) + \alpha_k\E_{v_k}[\nabla\obj(w_k) - \tilde{\phi}_{k,1} | v_{[k-1]}] \notag \\
            &\quad + (1-\alpha_k)(\nabla\obj(w_k) - \nabla\obj_{\delta}(w_k)),
    \end{align}
    where (s.1) follows by adding and subtracting $-\alpha_k\nabla\obj(w_k)$ and using \eqref{eq:unbiased_grad_est} in \cref{lem:grad_est_property}, and (s.2) relies on adding and subtracting $-(1-\alpha_k)\nabla\obj(w_k)$. By incorporating \eqref{eq:cond_expt_grad_est} into \eqref{eq:expectation_mid_term} and utilizing the tower rule, we obtain
    \begin{align}
        &\E_{v_{[k]}}[-\nabla \obj(w_k)^{\top}\tilde{\phi}_k] \notag \\
            &\!\leq \!-\!\E_{v_{[k]}}[\|\nabla\obj(w_k)\|^2] \!+\! \alpha_k\underbrace{\E_{v_{[k]}}[\nabla\obj(w_k)^{\top}(\nabla\obj(w_k) \!-\! \tilde{\phi}_{k,1})]}_{\numcircled{1}} \notag \\
            &\quad \!+\! (1\!-\!\alpha_k)\underbrace{\E_{v_{[k]}}[\nabla\obj(w_k)^{\top}(\nabla\obj(w_k) \!-\! \nabla\obj_{\delta}(w_k))]}_{\numcircled{2}}. \label{eq:cross_expt_midway}  
    \end{align}

    For term \numcircled{1} in \eqref{eq:cross_expt_midway}, we have
    \begin{align}
        \numcircled{1}
            &\stackrel{\text{(s.1)}}{\leq} \frac{1}{2} \E_{v_{[k]}}\left[\|\nabla\obj(w_k)\|^2 + \|\nabla\obj(w_k) \!-\! \nabla\obj(u_k) \!+\! \epsilon_k\|^2\right] \notag \\
            &\leq \frac{1}{2}\E_{v_{[k]}}[\|\nabla\obj(w_k)\|^2] \notag \\
            &\quad + \E_{v_{[k]}}[\|\nabla\obj(w_k) - \nabla\obj(u_k)\|^2] + \E_{v_{[k]}}[\|\epsilon_k\|^2] \notag \\
            &\stackrel{\text{(s.2)}}{\leq} \frac{1}{2}\E_{v_{[k]}}[\|\nabla\obj(w_k)\|^2] + L^2\delta^2 + \E_{v_{[k]}}[\|\epsilon_k\|^2], \label{eq:cross_expt_1}
    \end{align}
    where (s.1) leverages the inequality $\forall a,b \in \mathbb{R}^p,~a^\top b \leq 1/2\left(\|a\|^2 + \|b\|^2\right)$ and the definition of $\epsilon_k$, see \eqref{eq:grad_est_error}; 
    (s.2) relies on the $L$-smoothness of $\obj$ and the fact that $\|v_k\|=1, \forall v_k \sim U(\mathbb{S}_{p-1})$.

    The upper bound on term \numcircled{2} in \eqref{eq:cross_expt_midway} is
    \begin{align}
        \numcircled{2} &\leq \frac{1}{2} \Big(\E_{v_{[k]}}[\|\nabla\obj(w_k)\|^2] \!+\! \E_{v_{[k]}}[\|\nabla\obj(w_k) \!-\! \nabla\obj_{\delta}(w_k)\|^2]\Big) \notag \\
            &\leq \frac{1}{2}\E_{v_{[k]}}[\|\nabla\obj(w_k)\|^2] + \frac{L^2p^2\delta^2}{8}, \label{eq:cross_expt_2}
    \end{align}
    where the last inequality follows from \eqref{eq:difference_grad}. 

    We incorporate \eqref{eq:cross_expt_1} and \eqref{eq:cross_expt_2} into \eqref{eq:cross_expt_midway}. Hence, \eqref{eq:cross_expt_overall} holds.
    \end{proof}

    We provide an upper bound on the partial sum of a nonnegative sequence based on its recursive inequality.
    \begin{lemma}\label{lem:ubd_dec_seq_sum}
        If a nonnegative sequence $(g_k)_{k\in \mathbb{N}}$ satisfies $g_k \leq m_k g_{k-1} + \zeta_k, \forall k \in \mathbb{N}_{+}$, where $m_k \in (0,\bar{m}]$ and $\bar{m} \in (0,1)$, then,
        \begin{equation}\label{eq:decreasing_seq_sum_bd}
            \sum_{k=0}^{T} g_k \leq \frac{1}{1-\bar{m}} \left(g_0 + \sum_{k=1}^{T} \zeta_k\right).
        \end{equation}
    \end{lemma}
    \begin{proof}
        From the condition $m_k \in (0,\bar{m}]$ and the non-negativity of $(g_k)_{k\in \mathbb{N}}$, we know $g_k \leq \bar{m} g_{k-1} + \zeta_k, \forall k \in \mathbb{N}_{+}$.
        We sum over both sides of this inequality for $k=1,\ldots,T$ and obtain
        \begin{equation*}
            G_T - g_0 \leq \bar{m} (G_T - g_T) + \sum_{k=1}^{T} \zeta_k \leq \bar{m} G_T + \sum_{k=1}^{T} \zeta_k,
        \end{equation*}
        where $G_T \triangleq \sum_{k=0}^{T} g_k$. Hence, $(1-\bar{m}) G_T \leq g_0 + \sum_{k=1}^{T} \zeta_k$.
        Since $\bar{m}\in (0,1)$, the inequality \eqref{eq:decreasing_seq_sum_bd} holds.
    \end{proof}

    We quantify the order of the scaled cumulative error of inexact gradients $\sum_{k=0}^{T-1} \alpha_k \E_{v_{[T]}}[\|\epsilon_k\|^2]$ incurred by our gray-box controller~\eqref{eq:hybrid_controller}. 
    Here $T\in \mathbb{N}_{+}$ is the number of iterations, and $\E_{v_{[T]}}[\cdot]$ denotes the expectation with respect to the collection of i.i.d.~samples $v_{[T]} \triangleq \{v_0,\ldots,v_{T-1}\}$. 
    \begin{lemma}\label{lem:error_accumulation}
        Let \cref{assump:objective_Lipschitz} hold. If the approximate sensitivity $\hat{H}_k$ satisfies \eqref{eq:bound_sensitivity_fast_accurate}, \eqref{eq:bound_sensitivity_constant}, or \eqref{eq:bound_sensitivity_decay}, then the gray-box controller \eqref{eq:hybrid_controller} with the combination coefficients given in \eqref{eq:adaptive_comb_coeff_case_accurate}, \eqref{eq:adaptive_comb_coeff_case_constant}, or \eqref{eq:adaptive_comb_coeff_case_decay}, respectively, ensures that the error defined in \eqref{eq:grad_est_error} satisfies
        \begin{equation*}
            \sum_{k=0}^{T-1} \alpha_k \E_{v_{[T]}}[\|\epsilon_k\|^2] \!=\!
            \begin{cases}
                \bigO{T^{\max(0, 1-2\theta)}} & \text{if \eqref{eq:bound_sensitivity_fast_accurate}, \eqref{eq:adaptive_comb_coeff_case_accurate} hold, $\theta \!\neq\! \frac{1}{2}$, } \\
                \bigO{\ln T} & \text{if \eqref{eq:bound_sensitivity_fast_accurate}, \eqref{eq:adaptive_comb_coeff_case_accurate} hold, $\theta \!=\! \frac{1}{2}$,} \\
                \bigO{pT^{\frac{2}{3}}} & \text{if \eqref{eq:bound_sensitivity_constant}, \eqref{eq:adaptive_comb_coeff_case_constant} hold,} \\
                \bigO{pT^{\frac{2}{3} - \theta}} & \text{if \eqref{eq:bound_sensitivity_decay}, \eqref{eq:adaptive_comb_coeff_case_decay} hold.}  
            \end{cases}
        \end{equation*}
    \end{lemma}
    \begin{proof}
        First, we address the case when $\hat{H}_k$ satisfies \eqref{eq:bound_sensitivity_decay} and $\alpha_k$ follows \eqref{eq:adaptive_comb_coeff_case_decay}. The squared norm of the error $\epsilon_k$ satisfies
        \begin{align}\label{eq:bound_square_error_grad}
            \|\epsilon_k\|^2 &\stackrel{\text{(s.1)}}{=} \|(H_k \!-\! \hat{H}_k) \nabla_y \Phi(u_k,y_k)\|^2 \notag \\
                &\stackrel{\text{(s.2)}}{\leq} \epsilon_{H,k}^2 \|\nabla_y \Phi(u_k,y_k)\|^2 \stackrel{\text{(s.3)}}{\leq} \epsilon_{H,k}^2 M_\Phi^2 \stackrel{\text{(s.4)}}{\leq} \frac{M_\Phi^2 (\bar{\epsilon}')^2}{(k+1)^{2\theta}},
        \end{align}
        where (s.1) uses \eqref{eq:grad_est_error}; (s.2) follows from the inequality $\forall A\in \mathbb{R}^{p\times q}, b\in \mathbb{R}^q, \|Ab\| \leq \|A\|\|b\|$; (s.3) holds because the property that $\Phi(u,y)$ is $M_\Phi$-Lipschitz in $y$ (see \cref{assump:objective_Lipschitz}) implies that $\|\nabla_y \Phi(u,y)\| \leq M_\Phi$; (s.4) is due to \eqref{eq:bound_sensitivity_decay}. 
        Therefore,
        \begin{align*}
            \sum_{k=0}^{T-1} &\alpha_k \E_{v_{[T]}}[\|\epsilon_k\|^2] \leq \sum_{k=0}^{T-1} \alpha_k M_\Phi^2 \epsilon_{H,k}^2 \leq \sum_{k=0}^{T-1} \alpha_k \frac{M_\Phi^2 (\bar{\epsilon}')^2}{(k+1)^{2\theta}} \\
            &\stackrel{\text{(s.1)}}{\leq} M_\Phi^2 (\bar{\epsilon}')^2 \bigg(1 + \int_{0}^{T-1} \frac{C p}{(x+1)^{\frac{1}{3} + \theta}} \ud x \bigg) \\
            &\leq \bigg(1 + \frac{3C p}{2-3\theta}(T^{\frac{2}{3}-\theta} - 1) \bigg) M_\Phi^2 (\bar{\epsilon}')^2,
        \end{align*}
        where (s.1) uses the following inequality
        \begin{equation}\label{eq:sum_integral_test}
            \sum_{k=0}^{T-1} \frac{\alpha_k}{(k\!+\!1)^{2\theta}} \leq 1 \!+\! \sum_{k=1}^{T-1} \frac{C p}{(k\!+\!1)^{\frac{1}{3} + \theta}} \leq 1 \!+\! \int_0^{T-1} \! \frac{C p}{(x\!+\!1)^{\frac{1}{3} + \theta}} \ud x.
        \end{equation}

        For the case when $\hat{H}_k$ and $\alpha_k$ satisfy \eqref{eq:bound_sensitivity_constant} and \eqref{eq:adaptive_comb_coeff_case_constant}, respectively, we follow a similar reasoning and obtain
        \begin{align*}
            \sum_{k=0}^{T-1} &\alpha_k \E_{v_{[T]}}[\|\epsilon_k\|^2] \leq \sum_{k=0}^{T-1} \alpha_k M_\Phi^2 \bar{\epsilon}^2 \leq M_\Phi^2 \bar{\epsilon}^2 \bigg(1 + \sum_{k=1}^{T-1} \frac{Cp}{(k+1)^{\frac{1}{3}}} \bigg) \\
                &\leq \left(1 + \frac{3}{2} Cp (T^{\frac{2}{3}} \!-\! 1)\right) M_{\Phi}^2 \bar{\epsilon}^2.
        \end{align*}

        Finally, when the approximate sensitivity satisfies \eqref{eq:bound_sensitivity_fast_accurate} and the coefficient follows \eqref{eq:adaptive_comb_coeff_case_accurate}, we have
        \begin{align*}
            \sum_{k=0}^{T-1} &\alpha_k \E_{v_{[T]}}[\|\epsilon_k\|^2] \leq \sum_{k=0}^{T-1} \E_{v_{[T]}}[\|\epsilon_k\|^2] \leq \sum_{k=0}^{T-1} \frac{M_\Phi^2 (\bar\epsilon')^2}{(k+1)^{2\theta}} \notag \\
            &\leq 
            \begin{cases}
                M_\Phi^2 (\bar\epsilon')^2 \left(\frac{T^{1-2\theta}-1}{1-2\theta} + 1\right), & \text{ if } \theta \geq \frac{1}{3}, \theta \neq \frac{1}{2}, \\
                M_\Phi^2 (\bar\epsilon')^2 (\ln T + 1), & \text{ if } \theta = \frac{1}{2}.
            \end{cases}
        \end{align*}
        
        Hence, the orders of complexity in \cref{lem:error_accumulation} are proved.
    \end{proof}


\subsection{Proof of Lemma~\ref{lem:sec_mom_ubd}}\label{appendix:sec_mom}
    The second moment of $\tilde{\phi}_k$ in \eqref{eq:hybrid_controller_gradient} satisfies
    \begin{align}
        &\E_{v_{[k]}}[\|\tilde{\phi}_k\|^2]  = \E_{v_{[k]}}[\|\alpha_k\tilde{\phi}_{k,1} + (1-\alpha_k)\tilde{\phi}_{k,2}\|^2] \notag \\
            &\quad \leq 2\alpha_k^2 \E_{v_{[k]}}[\|\tilde{\phi}_{k,1}\|^2] + 2(1-\alpha_k)^2 \E_{v_{[k]}}[\|\tilde{\phi}_{k,2}\|^2], \label{eq:sec_mom_decomp}
    \end{align}
    where the last inequality follows from the fact that $\|a+b\|^2 \leq 2\|a\|^2 + 2\|b\|^2, \forall a,b \in \mathbb{R}^p$. 
    The bound on $\E_{v_{[k]}}[\|\tilde{\phi}_{k,1}\|^2]$ is
    \begin{align}
        &\E_{v_{[k]}}[\|\tilde{\phi}_{k,1}\|^2] \stackrel{\text{(s.1)}}{=} \E_{v_{[k]}}[\|\nabla\obj(u_k) - \epsilon_k\|^2] \notag \\
            &\ \leq 2\E_{v_{[k]}}[\|\nabla\obj(u_k)\|^2] + 2\E_{v_{[k]}}[\|\epsilon_k\|^2] \notag \\
            &\ = 2\E_{v_{[k]}}[\|\nabla\obj(u_k) \!-\! \nabla\obj(w_k) \!+\! \nabla\obj(w_k)\|^2] + 2\E_{v_{[k]}}[\|\epsilon_k\|^2] \notag \\
            &\ \leq 4\E_{v_{[k]}}[\|\nabla\obj(u_k) - \nabla\obj(w_k)\|^2] \notag \\
            &\qquad + 4\E_{v_{[k]}}[\|\nabla\obj(w_k)\|^2] + 2\E_{v_{[k]}}[\|\epsilon_k\|^2] \notag \\
            &\ \stackrel{\text{(s.2)}}{\leq} 4L^2\delta^2 + 4\E_{v_{[k]}}[\|\nabla\obj(w_k)\|^2] + 2\epsilon_{H,k}^2 M_\Phi^2, \label{eq:sec_mom_ig}
    \end{align}
    where (s.1) follows from the definition \eqref{eq:grad_est_error}; 
    \rev{(s.2) leverages $\|\nabla\obj(u_k) - \nabla\obj(w_k)\| \leq L\|\delta v_k\| = L\delta$ thanks to the $L$-smoothness of $\obj$, as well as (s.3) in \eqref{eq:bound_square_error_grad}}. 
    For the term $\E_{v_{[k]}}[\|\tilde{\phi}_{k,2}\|^2]$, we have
    \begin{align}
        &\E_{v_{[k]}}[\|\tilde{\phi}_{k,2}\|^2] \notag \\
            &= \frac{p^2}{\delta^2} \E_{v_{[k]}}[\|v_k(\obj(w_k+\delta v_k) - \obj(w_{k-1}+\delta v_{k-1}))\|^2] \notag \\
            &= \frac{p^2}{\delta^2} \E_{v_{[k]}}[\|v_k(\obj(w_k+\delta v_k) - \obj(w_{k-1}+\delta v_k)) \notag \\
            &\hspace{5em} + v_k(\obj(w_{k-1}+\delta v_k) - \obj(w_{k-1}+\delta v_{k-1}))\|^2] \notag \\
            &\leq \frac{2p^2}{\delta^2} \underbrace{\E_{v_{[k]}}[\|v_k(\obj(w_k\!+\!\delta v_k) - \obj(w_{k-1}\!+\!\delta v_k))\|^2]}_{\numcircled{1}} \notag \\
            &\quad + \frac{2p^2}{\delta^2} \underbrace{\E_{v_{[k]}}[\|v_k(\obj(w_{k-1}\!+\!\delta v_k) - \obj(w_{k-1}\!+\!\delta v_{k-1}))\|^2]}_{\numcircled{2}}. \label{eq:sec_mom_est_midway}
    \end{align}
    The upper bound on term $\numcircled{1}$ is given by
    \begin{align}
        \numcircled{1} \stackrel{\text{(s.1)}}{\leq} \E_{v_{[k]}}[\|v_k\|^2 M^2 \|w_k\!-\!w_{k\!-\!1}\|^2] \stackrel{\text{(s.2)}}{=} M^2\eta^2 \E_{v_{[k]}}[\|\tilde{\phi}_{k-1}\|^2], \label{eq:sec_mom_est_1}
    \end{align}
    where (s.1) is due to the assumption that $\obj(u)$ is $M$-Lipschitz, and (s.2) uses the update \eqref{eq:hybrid_controller_GD_update} and the fact that $\|v_k\|=1, \forall v_k\sim U(\mathbb{S}_{p-1})$.
    The upper bound on term $\numcircled{2}$ is
    \begin{align}
        &\numcircled{2} \stackrel{\text{(s.1)}}{=} \E_{v_{[k]}}[(\obj(w_{k-1}+\delta v_k) - \obj(w_{k-1}+\delta v_{k-1}))^2] \notag \\
            &\stackrel{\text{(s.2)}}{\leq} 3\E_{v_{[k]}}[(\obj(w_{k-1}+\delta v_k) - \obj(w_{k-1}) - \delta\nabla\obj(w_{k-1})^{\top}v_k)^2] \notag \\
            &\hspace{1em} + 3\E_{v_{[k]}}[(\obj(w_{k\!-\!1}) \!-\! \obj(w_{k\!-\!1} \!+\! \delta v_{k\!-\!1}) \!+\! \delta\nabla\obj(w_{k\!-\!1})^{\top}v_{k\!-\!1})^2] \notag \\
            &\hspace{1em} + 3\E_{v_{[k]}}[(\nabla\obj(w_{k-1})^{\top}\delta(v_k-v_{k-1}))^2] \notag \\
            &\stackrel{\text{(s.3)}}{\leq} 3\cdot\frac{L^2}{4}\E_{v_{[k]}}[\|\delta v_k\|^4] + 3\cdot\frac{L^2}{4}\E_{v_{[k]}}[\|\delta v_{k-1}\|^4] \notag \\
            &\hspace{1.2em} + 3\delta^2\E_{v_{[k]}}[\|\nabla\obj(w_{k-1})\|^2 \|v_k-v_{k-1}\|^2] \notag \\
            &\stackrel{\text{(s.4)}}{\leq} \frac{3}{2}\delta^4L^2 + 3\delta^2\E_{v_{[k]}}[\|\nabla\obj(w_{k-1})\|^2] \cdot \E_{v_{[k]}}[\|v_k - v_{k-1}\|^2] \notag \\
            &\stackrel{\text{(s.5)}}{\leq} 3\delta^2\Big(\frac{1}{2}\delta^2L^2 + 4\E_{v_{[k]}}[\|\nabla\obj(w_{k-1})\|^2]\Big), \label{eq:sec_mom_est_2}
    \end{align}
    where (s.1) holds since $\|v_k\|=1, \forall v_k\sim U(\mathbb{S}_{p-1})$; (s.2) is obtained by transforming $\obj(w_{k-1}+\delta v_k) - \obj(w_{k-1}+\delta v_{k-1})$ to the sum of three terms and then using the inequality $(a+b+c)^2 \leq 3(a^2+b^2+c^2), \forall a,b,c\in \mathbb{R}$; (s.3) uses the property of the $L$-smooth function $\obj(u)$ (see \cite[Eq.~(6)]{nesterov2017random}) and the Cauchy-Schwarz inequality; (s.4) follows from the independence between $\nabla\obj(w_{k-1})$ and $v_k,v_{k-1}$; (s.5) relies on the following bound
    \begin{equation*}
        \E_{v_{[k]}}[\|v_k - v_{k-1}\|^2] \leq 2\E_{v_{[k]}}[\|v_k\|^2] + 2\E_{v_{[k]}}[\|v_{k-1}\|^2] = 4.
    \end{equation*}
    By incorporating \eqref{eq:sec_mom_est_1} and \eqref{eq:sec_mom_est_2} into \eqref{eq:sec_mom_est_midway}, we have
    \begin{align}
        \E_{v_{[k]}}[\|\tilde{\phi}_{k,2}\|^2] \leq &\frac{2}{\delta^2}p^2M^2\eta^2 \E_{v_{[k]}}[\|\tilde{\phi}_{k-1}\|^2] + 3L^2p^2\delta^2 \notag \\
            &+ 24p^2\E_{v_{[k]}}[\|\nabla\obj(w_{k-1})\|^2]. \label{eq:sec_mom_est}
    \end{align}

    We plug \eqref{eq:sec_mom_ig} and \eqref{eq:sec_mom_est} into \eqref{eq:sec_mom_decomp}. Then, \eqref{eq:sec_mom_overall} holds.

\subsection{Proof of Theorem~\ref{thm:optimality_general}}\label{appendix:optimality_general}
    \revise{First, we present the convergence analysis for the gray-box controller \eqref{eq:hybrid_controller} with the combination coefficient $\alpha_k$ given by \eqref{eq:adaptive_comb_coeff_case_decay} when the sensitivity satisfies \eqref{eq:bound_sensitivity_decay}. The analyses for the remaining cases follow similar reasoning and will be presented later.}   
    Since the objective function $\obj$ is $L$-smooth, we have 
    \begin{align}\label{eq:smooth_func_ubd}
        \obj(w_{k+1}) &\leq \obj(w_k) \!+\! \nabla \obj(w_k)^{\top}(w_{k+1}\!-\!w_{k}) \!+\! \frac{L}{2} \|w_{k+1}\!-\!w_{k}\|^2 \notag \\
            &= \obj(w_k) - \eta\nabla \obj(w_k)^{\top}\tilde{\phi}_k + \frac{L\eta^2}{2}\|\tilde{\phi}_k\|^2.
    \end{align}
    We take expectations of both sides of \eqref{eq:smooth_func_ubd} with respect to $v_{[T]}$, sum them up for $k=0,\ldots,T-1$, and obtain
    \begin{align}\label{eq:smooth_func_expt_ubd_midway}
        &\E_{v_{[T]}}[\obj(w_T)] \leq \obj(w_0)  \notag \\
            & \quad +\! \eta\underbrace{\sum\limits_{k=0}^{T-1} \E_{v_{[T]}}[-\!\nabla\obj(w_k)^{\top}\tilde{\phi}_k]}_{\numcircled{1}} + \frac{L\eta^2}{2}\underbrace{\sum\limits_{k=0}^{T-1} \E_{v_{[T]}}[\|\tilde{\phi}_k\|^2]}_{\numcircled{2}}.
    \end{align}
    To derive an upper bound on the cross term \numcircled{1} in \eqref{eq:smooth_func_expt_ubd_midway}, we refer to \eqref{eq:cross_expt_overall} in \cref{lem:cross_expt_ubd}, exploit $\|\epsilon_k\|^2 \leq M_\Phi^2 \epsilon_{H,k}^2$ as \eqref{eq:bound_square_error_grad}, and obtain
    \begin{align}\label{eq:sum_cross_term_ubd}
        \numcircled{1} \leq& -\frac{1}{2} \sum_{k=0}^{T-1} \E_{v_{[T]}}[\|\nabla \obj(w_k)\|^2] + \sum_{k=0}^{T-1} \alpha_k M_{\Phi}^2 \epsilon_{H,k}^2 \notag \\
            & + L^2\delta^2 \sum_{k=0}^{T-1} \left(\alpha_k + \frac{p^2}{8}(1-\alpha_k)\right).
    \end{align}
    We construct an upper bound on term \numcircled{2} in \eqref{eq:smooth_func_expt_ubd_midway} as follows. 
    \revise{Given the combination coefficient $\alpha_k \in (0,1]$ and the specified parametric conditions of $\eta$ and $\delta$, the coefficient of $\E_{v_{[k]}}[\|\tilde{\phi}_{k-1}\|^2]$ on the right-hand side of the recursive inequality \eqref{eq:sec_mom_overall} satisfies
    \begin{equation*}
        \frac{4(1-\alpha_k)^2}{\delta^2}p^2M^2\eta^2 \leq \frac{4p^2M^2\eta^2}{\delta^2} \triangleq c, ~~ c = \frac{1}{2T^{\frac{2}{3}-2\theta}} \in \Big(0,\frac{1}{2}\Big].
    \end{equation*}
    Further, we utilize the following intermediate upper bound
    \begin{align}\label{eq:inter_cumul_sec_mom}
        &\sum_{k=1}^{T-1} \! \big(8\alpha_k^2 \E_{v_{[T]}}[\|\nabla\obj(w_k)\|^2] \!+\! 48(1\!-\!\alpha_k)^2p^2 \E_{v_{[T]}}[\|\nabla\obj(w_{k\!-\!1})\|^2]\big) \notag \\
            &\quad \stackrel{\text{(s.1)}}{\leq} \sum_{k=0}^{T-1} (8+48(1\!-\!\alpha_{k+1})p^2) \E_{v_{[T]}}[\|\nabla\obj(w_k)\|^2] \notag \\
            &\quad \stackrel{\text{(s.2)}}{\leq} \sum_{k=0}^{T-1} (8+48(1\!-\!\alpha_T)p^2) \E_{v_{[T]}}[\|\nabla\obj(w_k)\|^2],
    \end{align}
    where (s.1) holds because $\alpha_k \in (0,1]$ and the second moments are nonnegative; (s.2) is due to the fact that $\alpha_k$ is non-increasing with respect to $k$.} 
    Based on \eqref{eq:sec_mom_overall} in \cref{lem:sec_mom_ubd}, \eqref{eq:decreasing_seq_sum_bd} in \cref{lem:ubd_dec_seq_sum}, the condition that $\alpha_k \in (0,1]$, and \eqref{eq:inter_cumul_sec_mom}, we obtain
    \begin{align}\label{eq:sum_sec_mom_ubd}
        \numcircled{2} \leq &\frac{1}{1\!-\!c} \! \bigg[\! \E[\|\tilde{\phi}_0\|^2] \!+\! \sum_{k=0}^{T-1}(8 \!+\!48(1\!-\!\alpha_T)p^2) \! \E[\|\nabla\obj(w_k)\|^2] \notag \\
            &\!+\! \sum_{k=1}^{T-1}\!\left(4\alpha_k^2 M_\Phi^2 \epsilon_{H,k}^2 \!+\! 8\alpha_k^2L^2\delta^2 \!+\! 6(1\!-\!\alpha_k)^2L^2p^2\delta^2\right)\!\!\bigg],
    \end{align}

    We incorporate \eqref{eq:sum_cross_term_ubd} and \eqref{eq:sum_sec_mom_ubd} into \eqref{eq:smooth_func_expt_ubd_midway}, rearrange terms, utilize $\E_{v_{[T]}}[\obj(w_T)] \geq \obj^*$, and arrive at
    \begin{align}\label{eq:convergence_measure_midway}
        &\frac{1}{2}\Big(1\!-\!\frac{L\eta (8\!+\!48(1\!-\!\alpha_T)p^2)}{1\!-\!c}\Big) \cdot \frac{1}{T}\sum_{k=0}^{T-1}\E_{v_{[T]}}[\|\nabla\obj(w_k)\|^2] \notag \\
            &\leq \frac{\big(\obj(w_0) \!-\! \obj^*\big)}{\eta T} \!\!+\!\! \frac{1}{T}\! \bigg[\sum_{k=0}^{T\!-\!1}\alpha_k M_{\Phi}^2 \epsilon_{H,k}^2 \!+\! L^2\delta^2 \sum_{k=0}^{T\!-\!1}\Big(\alpha_k \!\!+\!\! \frac{p^2(1\!-\!\alpha_k)}{8}\Big)\bigg] \notag \\
            &\quad + \frac{L\eta}{2T(1-c)} \bigg[\E_{v_{[T]}}[\|\tilde\phi_0\|^2] + \sum_{k=1}^{T-1} 4\alpha_k^2 M_\Phi^2 \epsilon_{H,k}^2 \notag \\
            &\hspace{7em} + 2L^2\delta^2 \sum_{k=1}^{T-1} \big(4\alpha_k^2 \!+\! 3(1\!-\!\alpha_k)^2p^2\big) \bigg].
    \end{align}
    Recall that $c \!\in\! (0,\!\frac{1}{2}]$. For a sufficiently large $T$ that ensures $\alpha_T = Cp/(T\!+\!1)^{\frac{1}{3}-\theta}$, we have
    \begin{align*}
        &{\textstyle \frac{L\eta \left(8\!+\!48(1\!-\!\alpha_T)p^2\right)}{1-c}} \stackrel{\text{(s.1)}}{\leq} 2L\eta \Big(8 \!+\! \frac{48p^2}{\alpha_T}\Big) \!<\! \frac{(T\!+\!1)^{\frac{1}{3}-\theta}}{2T^{\frac{2}{3}-\theta}} \!<\! \frac{1}{2}, \\
        &\frac{L\eta}{2T(1-c)} \leq \frac{1}{224\max(C,\frac{1}{C})p T^{\frac{5}{3}-\theta}},
    \end{align*}
    where (s.1) uses the inequality $1-\alpha_T \leq \frac{1}{\alpha_T}-1<\frac{1}{\alpha_T}$ for $\alpha_T \in (0,1]$. 
    We proceed to characterize some important terms in \eqref{eq:convergence_measure_midway}:
    \begin{align*}
        &L^2\delta^2 \sum_{k=0}^{T-1} \alpha_k \leq L^2\delta^2 \sum_{k=0}^{T-1} \frac{C p}{(k\!+\!1)^{\frac{1}{3}\!-\!\theta}} 
        \stackrel{\text{(s.1)}}{\leq} CL^2 \delta^2 p \bigg(\! 1 \!+\! \frac{3(T^{\frac{2}{3}\!+\!\theta} \!-\! 1)}{2\!+\!3\theta}\!\bigg),
        \\
        &L^2\delta^2 p^2 \sum_{k=0}^{T-1} (1\!-\!\alpha_k) \stackrel{\text{(s.2)}}{\leq} \sum_{k=0}^{T-1} \frac{L^2\delta^2 p^2}{\alpha_k} \leq \sum_{k=0}^{T-1} \frac{L^2\delta^2 p (k+1)^{\frac{1}{3}-\theta}}{C} \\
        &\hspace{8.5em} \leq \frac{L^2 \delta^2 p \left((T+1)^{\frac{4}{3}-\theta}-1\right)}{C (\frac{4}{3}-\theta)},
    \end{align*}
    where (s.1) is similar to \eqref{eq:sum_integral_test}, and (s.2) again leverages $1-\alpha_k < \frac{1}{\alpha_k}$ when $\alpha_k \in (0,1]$. 
    Moreover, for any $k\in \mathbb{N}$, $\alpha_k \in (0,1]$, and $p,T\in \mathbb{N}_+$, we have the bounds $\alpha_k^2 M_\Phi^2 \epsilon_{H,k}^2 \leq \alpha_k M_\Phi^2 \epsilon_{H,k}^2$ and $4\alpha_k^2 \!+\! 3(1\!-\!\alpha_k)^2p^2 \leq 4\alpha_k \!+\! 3(1\!-\!\alpha_k)p^2$.
    
    Consequently, we derive from \eqref{eq:convergence_measure_midway} and \cref{lem:error_accumulation} the following upper bound
    \begin{align}\label{eq:convergence_measure_bd_acc}
        &\frac{1}{T}\sum_{k=0}^{T-1}\E_{v_{[T]}}[\|\nabla\obj(w_k)\|^2] \leq \underbrace{\frac{4(\obj(w_0) \!-\! \obj^*)}{\eta T}}_{\sim \mathcal{O}\big(p/T^{\frac{1}{3}+\theta}\big)} \notag \\
        &+ \underbrace{\frac{4}{T} \sum_{k=0}^{T-1} \alpha_k M_\Phi^2 \epsilon_{H,k}^2}_{\sim \mathcal{O}\big(p/T^{\frac{1}{3}+\theta}\big)} + \underbrace{\frac{L^2 \delta^2}{2T} \sum_{k=0}^{T-1} \left(8\alpha_k \!+\! p^2(1\!-\!\alpha_k)\right)}_{\sim \mathcal{O}\big(p/T^{\frac{1}{3}+\theta}\big)} \notag \\
        &+\!\! \underbrace{\frac{4L\eta}{T} \! \bigg[\! \E[\|\tilde{\phi}_0\|^2] \!+\! 4\! \sum_{k=1}^{T-1}\! \alpha_k \! M_\Phi^2 \epsilon_{H,k}^2 \!+\! 2L^2\delta^2 \!\! \sum_{k=1}^{T-1} \!\! \left(4\alpha_k \!+\! 3(1\!-\!\alpha_k)p^2\right) \!\! \bigg]}_{\leq C_{M,\|\tilde{\phi}_{0}\|^2}/T}\!.
    \end{align}
    The convergence measure grows on the order of $\mathcal{O}\big(p/T^{\frac{1}{3}+\theta}\big)$.
    
    The analysis when sensitivities are affected by bounded errors \eqref{eq:bound_sensitivity_constant} and $\alpha_k$ satisfies \eqref{eq:adaptive_comb_coeff_case_constant} follows a similar template as above, although $\theta$ is replaced with $0$. Therefore, the overall convergence measure admits the same upper bound as \eqref{eq:convergence_measure_bd_acc}, although the corresponding order becomes $\mathcal{O}(p/T^\frac{1}{3})$.

    When sensitivities are highly accurate as per \eqref{eq:bound_sensitivity_fast_accurate}, the parametric conditions of $\eta$ and $\delta$ still ensure that
    \begin{equation*}
        \frac{4(1-\alpha_k)^2}{\delta^2}p^2M^2\eta^2 \leq \frac{4p^2M^2\eta^2}{\delta^2} = c \in \Big(0,\frac{1}{2}\Big].
    \end{equation*}
    Further, different from \eqref{eq:inter_cumul_sec_mom}, the intermediate upper bound becomes
    \begin{align*}
        &\sum_{k=1}^{T-1} \! \big(8\alpha_k^2 \E_{v_{[T]}}[\|\nabla\obj(w_k)\|^2] \!+\! 48(1\!-\!\alpha_k)^2p^2 \E_{v_{[T]}}[\|\nabla\obj(w_{k\!-\!1})\|^2]\big) \notag \\
            &\quad \leq \sum_{k=0}^{T-1} (8\alpha_k+48(1\!-\!\alpha_{k+1})p^2) \E_{v_{[T]}}[\|\nabla\obj(w_k)\|^2] \notag \\
            &\quad \stackrel{\text{(s.1)}}{\leq} \sum_{k=0}^{T-1} (8+48pC) \E_{v_{[T]}}[\|\nabla\obj(w_k)\|^2],
    \end{align*}
    where (s.1) holds because $\forall k\in \mathbb{N}, 8\alpha_k\!+\!48(1\!-\!\alpha_{k+1})p^2 \leq 8\!+\!48Cp/(k\!+\!2)^{\frac{1}{3}} \leq 8\!+\!48Cp$. It follows that $\frac{L\eta(8+48pC)}{1-c} \leq \frac{1}{2T^{\frac{1}{3}}} \leq \frac{1}{2}$.
    Hence, we obtain the same upper bound as \eqref{eq:convergence_measure}. In this scenario,
    \begin{align*}
        &\frac{1}{\eta T} = \bigO{\frac{p}{T^\frac{2}{3}}}, \\
        &\frac{1}{T} \sum_{k=0}^{T-1} \alpha_k M_\Phi^2 \epsilon_{H,k}^2 = 
        \begin{cases}
            \bigO{\frac{1}{T^{\min(1, 2\theta)}}} \!=\! o\Big(\frac{1}{T^{\frac{2}{3}}}\Big), & \text{if } \theta \!\geq\! \frac{1}{3}, \theta \!\neq\! \frac{1}{2}, \\
            \bigO{\frac{\ln T}{T}} \!=\! o\Big(\frac{1}{T^{\frac{2}{3}}}\Big), & \text{if } \theta = \frac{1}{2},
        \end{cases}
        \\
        &\frac{L^2\delta^2}{T} \sum_{k=0}^{T-1} (8\alpha_k \!+\! p^2(1\!-\!\alpha_k)) \leq 8L^2\delta^2 \!+\! \frac{C L^2\delta^2p}{T} \sum_{k=0}^{T-1} \frac{1}{(k\!+\!1)^\frac{1}{3}} \\
        &\hspace{3em} = \bigO{\max\Big(\frac{1}{T^\frac{2}{3}}, \frac{p}{T}\Big)},
    \end{align*}
    where $\theta \geq \frac{1}{3}$. The convergence measure is therefore of the order $\mathcal{O}\big(p/T^\frac{2}{3}\big)$. Hence, the overall complexity results \eqref{eq:convergence_measure_complexity} hold. 

\subsection{Proof of Theorem~\ref{thm:model_based_optimality}}\label{appendix:model_based_optimality}
    \rev{The deterministic model-based feedback optimization controller is a special case of \eqref{eq:hybrid_controller}, using only the inexact gradient $\tilde{\phi}_{k,1}$ (i.e., with combination coefficients $\alpha_k=1,\forall k \in \mathbb{N}$) and updating without exploration signals (i.e., with a smoothing parameter $\delta=0$). We provide a separate proof to obtain sharper results.}

    In this case, $\forall k \in \mathbb{N}, u_k = w_k, \tilde{\phi}_k = \tilde{\phi}_{k,1}$. Hence, we use the definition \eqref{eq:grad_est_error} of $\epsilon_k$ to obtain
    \begin{equation}\label{eq:mb_sec_mom_bd}
        \|\tilde{\phi}_k\|^2 \leq 2 \|\nabla \obj(w_k)\|^2 + 2\|\epsilon_k\|^2.
    \end{equation}
    Furthermore, similar to \cref{appendix:cross_expt}, we have
    \begin{align}\label{eq:mb_cross_term_bd}
        -\nabla \obj(w_k)^{\top}\tilde{\phi}_k &\stackrel{\text{(s.1)}}{=} -\|\nabla \obj(w_k)\|^2 + \nabla\obj(w_k)^{\top}\epsilon_k \notag \\
            &\leq -\frac{1}{2} \|\nabla \obj(w_k)\|^2 + \frac{1}{2} \|\epsilon_k\|^2,
    \end{align}
    where (s.1) uses the definition \eqref{eq:grad_est_error} of $\epsilon_k$. 
    In the deterministic case, we telescope \eqref{eq:smooth_func_ubd} for $k=0,\ldots,T-1$ and obtain an inequality similar to \eqref{eq:smooth_func_expt_ubd_midway}, albeit without expectation. Then, we incorporate \eqref{eq:mb_sec_mom_bd} and \eqref{eq:mb_cross_term_bd}, rearrange terms, use $\obj(w_T) \geq \obj^*$, and arrive at
    \begin{align*}
        \Big(\frac{1}{2} \!-\! L\eta\Big) &\frac{1}{T}\sum_{k=0}^{T-1}\|\nabla\obj(w_k)\|^2 \leq \underbrace{\frac{\big(\obj(w_0) \!-\! \obj^*\big)}{\eta T}}_{\sim \bigO{1/T}} \!+\! \frac{1\!+\!2L\eta}{2T} \sum_{k=0}^{T-1} \|\epsilon_k\|^2.
    \end{align*}
    The parametric condition of $\eta$ implies $\tfrac{1}{4} \leq \frac{1}{2} - L\eta < \tfrac{1}{2}$ and that $\tfrac{1+2L\eta}{1-2L\eta} \leq 3$. We further use $\|\epsilon_k\|^2 \leq M_\Phi^2 \epsilon_{H,k}^2$ in \eqref{eq:bound_square_error_grad} to obtain \eqref{eq:convergence_measure_model_based}.

    Moreover, when the approximate sensitivity $\hat{H}_k$ satisfies \eqref{eq:bound_sensitivity_fast_accurate}, 
    \begin{align*}
        \sum_{k=0}^{T-1} \|\epsilon_k\|^2 &\stackrel{\text{(s.1)}}{\leq} \sum_{k=0}^{T-1} \frac{\bar{\epsilon}'^2 M^2_\Phi}{(k+1)^{2\theta}} \stackrel{\text{(s.2)}}{\leq} 
        \begin{cases}
            \bar{\epsilon}'^2 M^2_\Phi \big(\frac{T^{1-2\theta} \!-\! 1}{1-2\theta} \!+\! 1 \big), \hspace{-1.5ex} &\text{if } \theta \neq \frac{1}{2}, \\
            \bar{\epsilon}'^2 M^2_\Phi (\ln T + 1), &\text{if } \theta = \frac{1}{2},
        \end{cases}
    \end{align*}
    where (s.1) is similar to \eqref{eq:bound_square_error_grad}, and (s.2) follows analogously as \eqref{eq:sum_integral_test}. Hence, the right-hand side of \eqref{eq:convergence_measure_model_based} scales as $\mathcal{O}(1/T^{\min(1,2\theta)})$ for $\theta \geq 1/3, \theta \neq 1/2$, and as $\bigO{\ln T / T}$ for $\theta = 1/2$. 
    If \eqref{eq:bound_sensitivity_constant} holds, then $\sum_{k=0}^{T-1} \|\epsilon_k\|^2 = \bar{\epsilon}^2 T$. If \eqref{eq:bound_sensitivity_decay} holds, then $\sum_{k=0}^{T-1} \|\epsilon_k\|^2 = \bigO{T^{1-2\theta}}$, where $\theta \in (0,1/3)$. Hence, the complexity results in \eqref{eq:convergence_measure_complexity_MB} hold.


\subsection{Proof of Theorem~\ref{thm:optimality_general_mf}}\label{appendix:model_free_optimality}
The convergence analysis herein resembles \cref{appendix:optimality_general}, although $\alpha_k=0, \forall k \in \mathbb{N}$. We obtain an inequality similar to \eqref{eq:convergence_measure_midway} as follows:
\begin{align*}
    \frac{1}{2} &\Big(1 - \frac{L\eta(8\!+\!48p^2)}{1-c}\Big) \cdot \frac{1}{T}\sum_{k=0}^{T-1} \E_{v_{[T]}}[\|\nabla \obj(w_k)\|^2] \\
    &\leq \underbrace{\frac{\big(\obj(w_0) \!-\! \obj^*\big)}{\eta T}}_{\sim \mathcal{O}\big(p^2/T^\frac{2}{3}\big)} + \underbrace{\frac{L^2p^2\delta^2}{8}}_{\sim \mathcal{O}\big(p^2/T^\frac{2}{3}\big)} \\
    &\quad + \underbrace{\frac{L\eta}{2T(1\!-\!c)} \Big[\E_{v_{[T]}}[\|\tilde\phi_0\|^2] + 6L^2\delta^2p^2T \Big]}_{\sim \bigO{1/T}}.
\end{align*}
The parametric conditions ensure $c\in (0,\frac{1}{2}]$, $1-\frac{L\eta}{(1-c)}(8\!+\!48p^2) \in [\frac{1}{2},1)$, and $\frac{L\eta}{2T(1-c)} \in \big(0,\frac{1}{224p^2T^\frac{4}{3}}\big]$. Therefore, \cref{thm:optimality_general_mf} is proved.

\subsection{Proof of Lemma~\ref{lem:hybrid_prod_mom_bd}}
    First, we know from the combination rule \eqref{eq:hybrid_gradient_TV} that
    \begin{align}\label{eq:cross_term_TV}
        &\E_{v_{[k]}}[\hat{\phi}_k^\top(w_k\!-\!u_k^*)] \notag \\
            &= \alpha_k\underbrace{\E_{v_{[k]}}[\hat{\phi}_{k,1}^\top(w_k\!-\!u_k^*)]}_{\numcircled{1}} + (1\!-\!\alpha_k)\underbrace{\E_{v_{[k]}}[\hat{\phi}_{k,2}^\top(w_k\!-\!u_k^*)]}_{\numcircled{2}}.
    \end{align}
    For term \numcircled{1} in \eqref{eq:cross_term_TV},
    \begin{align}\label{eq:cross_term_lb_MB}
        \numcircled{1} &= \E_{v_{[k]}}[\nabla \obj_k(w_k)^{\top}(w_k - u_k^*)] \notag \\
            &\hspace{2em} + \E_{v_{[k]}}[(\hat{\phi}_{k,1} - \nabla\obj_k(w_k))^{\top}(w_k-u_k^*)] \notag \\
            & \stackrel{\text{(s.1)}}{\geq} \E_{v_{[k]}}[\obj_k(w_k)] \!-\! \obj_k(u_k^*) \!-\! D \underbrace{\E_{v_{[k]}}[\|\hat{\phi}_{k,1} \!\!-\!\! \nabla\obj_k(w_k)\|]}_{\numcircled{3}},
    \end{align}
    where (s.1) utilizes the convexity of $\obj_k$, the Cauchy-Schwarz inequality, and the inequality $\|w_k - u_k^*\| \leq D$ (see \cref{assump:constraint_set}). Furthermore, term \numcircled{3} in \eqref{eq:cross_term_lb_MB} satisfies
    \begin{align}\label{eq:upper_bd_dist_grad_first_part}
        \numcircled{3} &\stackrel{\text{(s.1)}}{\leq} \E_{v_{[k]}}[\|\hat{\phi}_{k,1} \!-\! \nabla\obj_k(u_k)\|] + \E_{v_{[k]}}[\|\nabla\obj_k(u_k) \!-\! \nabla\obj_k(w_k)\|] \notag \\
            &\stackrel{\text{(s.2)}}{\leq} \E_{v_{[k]}}[\|\epsilon_k\|] + L_k \delta \stackrel{\text{(s.3)}}{\leq} M_{\Phi,k} \epsilon_{H,k} + L_k \delta,
    \end{align}
    where (s.1) is obtained by adding and subtracting $\nabla\obj_k(u_k)$ and using the triangle inequality; (s.2) utilizes the shorthand \eqref{eq:grad_est_error}, the assumption that $\obj_k(u)$ is $L_k$-smooth, and the rule \eqref{eq:hybrid_TV_noise}; (s.3) follows similarly as parts of \eqref{eq:bound_square_error_grad}, i.e.,
    \begin{equation}\label{eq:bound_error_grad_TV}
        \|\epsilon_k\| = \|(H_k' - \hat{H}_k) \nabla_y \Phi_k(u_k,y_k)\| \leq \epsilon_{H,k} M_{\Phi,k}.
    \end{equation}
    For term \numcircled{2} in \eqref{eq:cross_term_TV}, we have
    \begin{align}\label{eq:cross_term_lb_MF}
        \numcircled{2} &\stackrel{\text{(s.1)}}{=} \E_{v_{[k-1]}}\big[\E_{v_k}[\hat{\phi}_{k,2}^{\top}(w_k-u_k^*)|v_{[k-1]}]\big] \notag \\
            &\stackrel{\text{(s.2)}}{=} \E_{v_{[k-1]}}\big[\E_{v_k}[\hat{\phi}_{k,2}|v_{[k-1]}]^{\top}(w_k-u_k^*)\big] \notag \\
            &\stackrel{\text{(s.3)}}{=} \E_{v_{[k-1]}}\big[\nabla \obj_{k,\delta}(w_k)^{\top}(w_k-u_k^*)\big] \notag \\
            &\stackrel{\text{(s.4)}}{\geq} \E_{v_{[k]}}\big[\obj_{k,\delta}(w_k) - \obj_{k,\delta}(u_k^*)\big] \notag \\
            &= \E_{v_{[k]}}\big[\obj_{k,\delta}(w_k) - \obj_k(w_k)\big] + \E_{v_{[k]}}\big[\obj_k(w_k) - \obj_k(u_k^*)\big] \notag \\
            &\quad + \E_{v_{[k]}}\big[\obj_k(u_k^*) - \obj_{k,\delta}(u_k^*)\big] \notag \\
            &\stackrel{\text{(s.5)}}{\geq} \E_{v_{[k]}}[\obj_k(w_k)] - \obj_k(u_k^*) - L\delta^2,
    \end{align}
    where (s.1) uses the tower rule; (s.2) holds since $w_k-u_k^*$ is measurable with respect to $v_{[k-1]}$; (s.3) follows from \eqref{eq:unbiased_grad_est} and the independence of $\Phi_{k-1}(u_{k-1},y_{k-1})$ and $v_k$; (s.4) uses the convexity of $\obj_{k,\delta}$ and the independence of $\obj_{k,\delta}(w_k)$ and $v_k$; (s.5) follows from \eqref{eq:difference_func_val}.
    We combine the bounds \eqref{eq:cross_term_lb_MB}, \eqref{eq:upper_bd_dist_grad_first_part}, and \eqref{eq:cross_term_lb_MF}, merge them into \eqref{eq:cross_term_TV}, and arrive at \eqref{eq:hybrid_prod_lbd}.

    We proceed to analyze the second moment $\E_{v_{[k]}}[\|\hat{\phi}_k\|^2]$. Based on the combination rule \eqref{eq:hybrid_gradient_TV}, we know
    \begin{equation}\label{eq:sec_mom_mid_bd}
        \E_{v_{[k]}}\![\|\hat{\phi}_k\|^2] \leq 2\alpha_k^2 \!\underbrace{\E_{v_{[k]}}\![\|\hat{\phi}_{k,1}\|^2]}_{\numcircled{1}} + 2(1\!-\!\alpha_k)^2 \!\underbrace{\E_{v_{[k]}}\![\|\hat{\phi}_{k,2}\|^2]}_{\numcircled{2}}.
    \end{equation}
    Term \numcircled{1} in \eqref{eq:sec_mom_mid_bd} satisfies $\numcircled{1} \leq 2(M_k^2 \!+\! M_{\Phi,k}^2 \epsilon_{H,k}^2)$, where we use \eqref{eq:bound_error_grad_TV} and $\|\hat{\phi}_{k,1}\|^2 = \|\nabla \obj_k(u_k) - \epsilon_k\|^2 \leq 2M_k^2 + 2\|\epsilon_k\|^2$, because the $M_k$-Lipschitz continuity of $\obj_k$ implies that $\forall u_k \in \mathcal{U}, \|\nabla \obj_k(u_k)\| \leq M_k$. Term \numcircled{2} in \eqref{eq:sec_mom_mid_bd} admits the following bound
    \begin{align*}
        &\E_{v_{[k]}}[\|\hat{\phi}_{k,2}\|^2] \notag \\
          &\stackrel{\text{(s.1)}}{=} \frac{p^2}{\delta^2}\E_{v_{[k]}}[|\obj_k(w_k\!+\!\delta v_k) - \obj_{k-1}(w_{k-1}\!+\!\delta v_{k-1})|^2] \notag \\
          &\leq \frac{2p^2}{\delta^2} \Big(\E_{v_{[k]}}\big[|\obj_k(w_k\!+\!\delta v_k)|^2 + |\obj_{k-1}(w_{k-1}\!+\!\delta v_{k-1})|^2\big] \Big) \notag \\
          &\stackrel{\text{(s.2)}}{\leq} \frac{4p^2G^2}{\delta^2},
    \end{align*}
    where (s.1) holds since $\forall v_k \sim U(\mathbb{S}_{p-1}),\|v_k\| = 1$, and (s.2) uses the boundedness of $\obj_k$, i.e., $\forall u \in \mathcal{U}_{\tau}, k \in \mathbb{N}, |\obj_k(u)| \leq G$. We incorporate these bounds into \eqref{eq:sec_mom_mid_bd} and obtain \eqref{eq:hybrid_mom_ubd}.

\subsection{Proof of Theorem~\ref{thm:track_dynamic_regret}}\label{appendix:track_dynamic_regret}
    Because the optimal point $u_k^*$ lies in $\mathcal{U}$, we know from \eqref{eq:hybrid_controller_TV_GD_update} and the Pythagorean theorem (see \cite[Theorem~2.1]{hazan2022introduction}) that
    \begin{align*}
        \|w_{k+1} - u_k^*\|^2 &\leq \|w_k - \eta\hat{\phi}_k - u_k^*\|^2 \\
            &= \|w_k - u_k^*\|^2 - 2\eta \hat{\phi}_k^{\top}(w_k - u_k^*) + \eta^2 \|\hat{\phi}_k\|^2.
    \end{align*}
    We rearrange terms and obtain
    \begin{equation}\label{eq:cross_term_TV_ub}
        \hat{\phi}_k^\top(w_k\!-\!u_k^*) \leq \frac{\eta}{2} \|\hat{\phi}_k\|^2 + \frac{\|w_k\!-\!u_k^*\|^2\!-\!\|w_{k+1}\!-\!u_k^*\|^2}{2\eta}.
    \end{equation}
    
    \revise{We combine the lower bound \eqref{eq:hybrid_prod_lbd} in \cref{lem:hybrid_prod_mom_bd} with \eqref{eq:cross_term_TV_ub} and telescope the inequality to obtain}
    \begin{align}\label{eq:dynamic_regret_intermediate}
        &\sum_{k=1}^{T} \left(\E_{v_{[k]}}[\obj_k(w_k)] - \obj_k(u_k^*)\right) \leq \frac{\eta}{2} \underbrace{\sum_{k=1}^{T} \E_{v_{[k]}}[\|\hat{\phi}_k\|^2]}_{\numcircled{1}} \notag \\
            &\quad + \frac{1}{2\eta} \underbrace{\sum_{k=1}^{T} \left(\E_{v_{[k]}}[\|w_k\!-\!u_k^*\|^2] - \E_{v_{[k]}}[\!\|w_{k+1}\!-\!u_k^*\|^2]\right)}_{\numcircled{2}} \notag \\
            &\quad + D\sum_{k=1}^{T} \alpha_k (M_{\Phi,k} \epsilon_{H,k} + L_k \delta) + \sum_{k=1}^{T} (1-\alpha_k)L_k \delta^2.
    \end{align}
    \revise{
    We know from \eqref{eq:hybrid_mom_ubd} in \cref{lem:hybrid_prod_mom_bd} that an upper bound on term \numcircled{1} in \eqref{eq:dynamic_regret_intermediate} is given by
    \begin{equation}\label{eq:dynamic_regret_sum_moment}
        \numcircled{1} \leq \sum_{k=1}^{T} \Big(4\alpha_k^2(M_k^2 \!+\! M_{\Phi,k}^2 \epsilon_{H,k}^2) \!+\! 8(1\!-\!\alpha_k)^2 \frac{p^2G^2}{\delta^2}\Big).
    \end{equation}
    }
    Furthermore, term \numcircled{2} in \eqref{eq:dynamic_regret_intermediate} satisfies
    \begin{align}\label{eq:dynamic_regret_sum_distance}
        \numcircled{2} &\leq \E_{v_{[k]}}[\|w_1\|^2] - 2\E_{v_{[k]}}[w_1^{\top}u_1^*] + 2\E_{v_{[k]}}[w_{T+1}^{\top}u_T^*] \notag \\
            &\quad + 2\sum_{k=1}^{T-1} \E_{v_{[k]}}[w_{k+1}^{\top}(u_k^* - u_{k+1}^*)] \notag \\
            &\stackrel{\text{(s.1)}}{\leq} 5\bar{D}^2 + 2\bar{D} \sum_{k=1}^{T-1} \|u_k^* - u_{k+1}^*\|,
    \end{align}
    where (s.1) uses the Cauchy-Schwarz inequality and the fact that $\forall u \in \mathcal{U}, \|u\| \leq \bar{D}$, see also the discussion below \cref{assump:objective_uniform_bound}. 
    \revise{We incorporate \eqref{eq:dynamic_regret_sum_moment} and \eqref{eq:dynamic_regret_sum_distance} into \eqref{eq:dynamic_regret_intermediate} and analyze the orders of the dynamic regret $\reg$ in different scenarios.}

    Let us consider the case when the asymptotically accurate sensitivity satisfies \eqref{eq:bound_sensitivity_decay}. In this scenario, based on the inequality $1-\alpha_k \leq \frac{1}{\alpha_k} - 1 < \frac{1}{\alpha_k}, \forall \alpha_k \in (0,1]$, the parametric conditions of $\eta$ and $\delta$, and the rule \eqref{eq:bound_sens_decay_TV}, we derive the following upper bounds
    \begin{align*}
        &4\eta \sum_{k=1}^{T}(1\!-\!\alpha_k)^2 \frac{p^2G^2}{\delta^2} \leq \frac{4\eta p^2 G^2}{\delta^2} \sum_{k=1}^{T} \frac{1}{\alpha_k^2} = \bigO{\sqrt{p}T^{\kappa}}, \\
        &\sum_{k=1}^{T} (1\!-\!\alpha_k)L_k \delta^2 \leq (\max_k L_k) \delta^2 \sum_{k=1}^{T} \frac{1}{\alpha_k} = \bigO{\sqrt{p}T^{\kappa}},
    \end{align*}
    where $\kappa=\max(\frac{5}{6}-\frac{\theta}{2}, \frac{3}{4})$. Furthermore, other related terms satisfy
    \begin{align*}
        &\eta \sum_{k=1}^{T}\alpha_k^2 M_{\Phi,k}^2\epsilon_{H,k}^2 \leq \eta \sum_{k=1}^{T} M_{\Phi,k}^2\epsilon_{H,k}^2 = \mathcal{O}\Big({p^{-\frac{1}{2}}T^{\min(\frac{1}{6}\!-\!\frac{3}{2}\theta, \frac{1}{4}\!-\!2\theta )}} \Big), \\
        &\sum_{k=1}^{T} \alpha_k  \epsilon_{H,k} \leq \sum_{k=1}^{T} \frac{C  p^\frac{1}{6}\bar{\epsilon}'}{(k\!+\!1)^{\max(\frac{1}{6}\!+\!\frac{\theta}{2}, \frac{1}{12}\!+\!\theta)}} \!=\! \bigO{p^\frac{1}{6}T^{\min(\frac{5}{6}\!-\!\frac{\theta}{2}, \frac{11}{12}\!-\!\theta)}}, \\
        &\sum_{k=1}^{T} \alpha_k L_k \delta = \bigO{\sqrt{p}T^{\min(\frac{2}{3}+\frac{\theta}{2}, \frac{3}{4})}}.
    \end{align*}
    All these terms are asymptotically smaller than $\bigO{\sqrt{p}T^\kappa}$ given $\theta \in (0,\frac{1}{3})$, where $\kappa=\max(\frac{5}{6}-\frac{\theta}{2}, \frac{3}{4})$.
    By incorporating \eqref{eq:dynamic_regret_sum_moment} and \eqref{eq:dynamic_regret_sum_distance} into \eqref{eq:dynamic_regret_intermediate} and invoking the above complexity results, we have the following upper bound on the dynamic regret
    \begin{align}\label{eq:dyn_reg_bd_acc}
        &\reg \!\leq\! \underbrace{\frac{1}{2\eta}(2\bar{D}C_T \!+\! 5\bar{D}^2)}_{\sim \bigO{\sqrt{p}T^\kappa (C_T + 1)}} \!+\! \underbrace{4\eta \sum_{k=1}^{T} \! (1\!-\!\alpha_k)^2 \frac{p^2G^2}{\delta^2} \!+\! \sum_{k=1}^{T} \! (1\!-\!\alpha_k)L_k \delta^2}_{\sim \bigO{\sqrt{p}T^{\kappa}}} \notag \\
        &+ \underbrace{2\eta\sum_{k=1}^{T} \alpha_k^2(M_k^2 \!+\! M_{\Phi,k}^2\epsilon_{H,k}^2) \!+\! D\sum_{k=1}^{T} \alpha_k M_{\Phi,k}\epsilon_{H,k} \!+\! D \sum_{k=1}^{T} \alpha_k L_k \delta}_{\sim o\left(\sqrt{p}T^{\kappa}\right)}.
    \end{align}
    The order of $\reg$ is therefore $\bigO{\sqrt{p}T^\kappa (C_T\!+\!1)}$.

    The case when the sensitivity suffers from bounded errors \eqref{eq:bound_sensitivity_constant} corresponds to setting $\theta=0$ in the above design and analysis. The upper bound \eqref{eq:dyn_reg_bd_acc} still holds. The resulting $\reg$ grows on the order of $\mathcal{O}\big(\sqrt{p}T^{\frac{5}{6}} (C_T\!+\!1)\big)$.

    We proceed to address the case when the sensitivity is highly accurate, i.e., \eqref{eq:bound_sensitivity_fast_accurate} holds. Let $\beta_k = \min \Big\{\frac{C}{p^{\frac{1}{6}}\sqrt{k+1}}, 1\Big\}$. Then, $\alpha_k=1-\beta_k, \forall k\in \mathbb{N}$. Analogously, we derive from \eqref{eq:dynamic_regret_intermediate}--\eqref{eq:dynamic_regret_sum_distance} and the parametric conditions in this case that the dynamic regret bound is
    \begin{align*}
        &\reg \leq \underbrace{\frac{1}{2\eta}(5\bar{D}^2+2\bar{D}C_T)}_{\sim \bigO{\sqrt{pT}(C_T\!+\!1)}} + \underbrace{2\eta\sum_{k=1}^{T} (M_k^2 + M_{\Phi,k}^2 \epsilon_{H,k}^2)}_{\sim \bigO{\sqrt{T}/\sqrt{p}}} \\
        &\quad +\! \underbrace{4\eta \frac{p^2G^2}{\delta^2} \sum_{k=1}^{T} \beta_k^2}_{\sim \bigO{\sqrt{pT}}} \!+\! \underbrace{D\sum_{k=1}^{T} L_k\delta}_{\sim \mathcal{O}\big(p^\frac{1}{3}\sqrt{T}\big)} \!+\! \underbrace{\sum_{k=1}^{T} \beta_k L_k \delta^2}_{\sim \bigO{\sqrt{p}/\sqrt{T}}} + D\sum_{k=1}^{T} M_{\Phi,k}\epsilon_{H,k}.
    \end{align*}
    The last term $D\sum_{k=1}^{T} M_{\Phi,k}\epsilon_{H,k}$ is of order $\bigO{T^{1-\theta}}$ when $\theta \in [\frac{1}{3},1)$, $\bigO{\ln T}$ when $\theta=1$, and $\bigO{1}$ when $\theta>1$. 
    In this scenario, $\reg$ scales as $\bigO{\sqrt{pT}(C_T\!+\!1) + T^{\max(\frac{1}{2}, 1-\theta)}}$. Hence, the overall order \eqref{eq:dynamic_regret_complexity} of the dynamic regret is proved.

\subsection{Proof of Theorem~\ref{thm:track_dynamic_regret_mf}}\label{appendix:track_dynamic_regret_mf}
Analogous to the analysis related to \eqref{eq:dynamic_regret_intermediate} in \cref{appendix:track_dynamic_regret}, the model-based controller (i.e., \eqref{eq:hybrid_controller_TV} with $\alpha_k=1,\forall k\in \mathbb{N}$ and $\delta=0$) incurs the following dynamic regret bound
\begin{align}\label{eq:dynamic_regret_mb_intermediate}
        \reg \leq \eta\sum_{k=1}^{T} \! (M_k^2 \!+\! M_{\Phi,k}^2 \epsilon_{H,k}^2) \!+\! \frac{2\bar{D}C_T \!+\! 5\bar{D}^2}{2\eta} \!+\! D\sum_{k=1}^{T} \! M_{\Phi,k}\epsilon_{H,k}.
    \end{align}
The parametric conditions ensure that the first two terms on the right-hand side of \eqref{eq:dynamic_regret_mb_intermediate} are of the order $\mathcal{O}(\sqrt{T}(C_T+1))$. The last term on the right-hand side of \eqref{eq:dynamic_regret_mb_intermediate} grows on the order of
\begin{align*}
    \sum_{k=1}^{T} M_{\Phi,k}\epsilon_{H,k} =
    \begin{cases}
        \bigO{T^{\max(0, 1-\theta)}} & \text{if \eqref{eq:bound_sensitivity_fast_accurate} holds, } \theta \!\geq\! \frac{1}{3}, \theta \!\neq\! 1, \\
        \bigO{\ln T} & \text{if \eqref{eq:bound_sensitivity_fast_accurate} holds, } \theta = 1, \\
        \bigO{T} & \text{if \eqref{eq:bound_sensitivity_constant} holds,} \\
        \bigO{T^{1-\theta}} & \text{if \eqref{eq:bound_sensitivity_decay} holds,}
    \end{cases}
\end{align*}
Hence, the complexity results in \eqref{eq:dynamic_regret_mb} are proved.

Similarly, the dynamic regret bound corresponding to the model-free controller (i.e., \eqref{eq:hybrid_controller_TV} with $\alpha_k=0,\forall k\in \mathbb{N}$) is
\begin{align}\label{eq:dynamic_regret_mf_intermediate}
    \reg \leq \frac{2G^2 \eta T p^2}{\delta^2} + \frac{1}{2\eta}(2\bar{D}C_T + 5\bar{D}^2) + LT\delta^2.
\end{align}
The specified parametric conditions ensure that the right-hand side of \eqref{eq:dynamic_regret_mf_intermediate} is of the order $\mathcal{O}\big(p^\frac{2}{3}T^\frac{2}{3}(C_T+1)\big)$.

\begin{table*}[!tb]
\centering
\renewcommand \arraystretch{0.5}
\caption{Orders of dynamic regret for problem \eqref{eq:opt_TV}}
\label{table:complexity_TV}
\begin{threeparttable}
    \begin{tabularx}{.95\linewidth}{*{4}{Y}}
        \toprule
        \makecell{{\bfseries Controllers}} & \makecell{Highly accurate \eqref{eq:bound_sensitivity_fast_accurate} \\ $\theta \geq \frac{1}{3}$} & \makecell{Asymptotically accurate \eqref{eq:bound_sensitivity_decay} \\ $\theta \in (0, \frac{1}{3})$} & Bounded errors \eqref{eq:bound_sensitivity_constant} \\
        \midrule
        \makecell{model-based \eqref{eq:hybrid_controller_TV} \\ $\alpha_k\!=\!1,\delta\!=\!0$} & $\bigO{\sqrt{T}C_T\!+\!T^{\max(\frac{1}{2},1-\theta)}}$ & $\bigO{\sqrt{T}(C_T \!+\! 1) \!+\! T^{1-\theta}}$ & $\bigO{\sqrt{T}(C_T \!+\! 1) \!+\! T}$ \\
        \midrule
        \makecell{model-free \eqref{eq:hybrid_controller_TV} \\ with $\alpha_k=0$} & $\bigO{p^{\frac{2}{3}}T^{\frac{2}{3}}(C_T \!+\! 1)}$ & $\bigO{p^{\frac{2}{3}}T^{\frac{2}{3}}(C_T \!+\! 1)}$ & $\bigO{p^{\frac{2}{3}}T^{\frac{2}{3}}(C_T \!+\! 1)}$ \\
        \midrule
        gray-box \eqref{eq:hybrid_controller_TV} & $\bigO{\!\sqrt{pT}(C_T\!+\!1) \!+\! T^{\max(\frac{1}{2}, 1\!-\!\theta)}\!}$ & $\bigO{\sqrt{p} T^{\max(\frac{5}{6}-\frac{\theta}{2},\frac{3}{4})}(C_T \!+\! 1)}$ & $\bigO{\sqrt{p} T^{\frac{5}{6}}(C_T \!+\! 1)}$ \\
        \bottomrule
    \end{tabularx}
\end{threeparttable}
\end{table*}

We summarize the dynamic regret characterizations of different controllers when applied to problem~\eqref{eq:opt_TV} in \cref{table:complexity_TV}.

\ifshowtracking
\subsection{Proof of Theorem~\ref{thm:track_performance}}\label{appendix:track_performance}
    The recursive relation of the tracking error is
    \begin{align}\label{eq:upper_bd_track_err}
        \|w_{k+1} &- u_{k+1}^*\| \stackrel{\text{(s.1)}}{=} \|\proj(w_k - \eta \hat{\phi}_k) - u_{k+1}^*\| \notag \\
            &\stackrel{\text{(s.2)}}{\leq} \|\proj(w_k - \eta \hat{\phi}_k) - u_k^*\| + \|u_{k+1}^* - u_k^*\|,
    \end{align}
    where (s.1) uses \eqref{eq:hybrid_controller_TV_GD_update}, and (s.2) follows by adding and subtracting $u_k^*$ and using the triangle inequality. For the first term in \eqref{eq:upper_bd_track_err}, we have the following upper bound
    \begin{align}\label{eq:upper_bd_track_err_first_term} 
        \|&\proj(w_k - \eta\hat{\phi}_k) - u_k^*\| \notag \\
            &\stackrel{\text{(s.1)}}{=} \|\proj(w_k - \eta\hat{\phi}_k) - \proj(u_k^* - \eta\nabla \obj_k(u_k^*))\| \notag \\
            &\stackrel{\text{(s.2)}}{\leq} \|w_k - \eta\hat{\phi}_k - (u_k^* - \eta\nabla \obj_k(u_k^*))\| \notag \\
            &\stackrel{\text{(s.3)}}{\leq} \|w_k - \eta \nabla\obj_k(w_k) - (u_k^* - \eta\nabla \obj_k(u_k^*))\| \notag \\
            &\quad + \eta\|\hat{\phi}_k - \nabla\obj_k(w_k)\| \notag \\
            &\stackrel{\text{(s.4)}}{\leq} c_k\|w_k - u_k^*\| + \eta\|\hat{\phi}_k - \nabla\obj_k(w_k)\|,
    \end{align}
    where $c_k \triangleq \max\{|1-\eta\mu_k|,|1-\eta L_k|\}$.
    In \eqref{eq:upper_bd_track_err_first_term}, (s.1) holds because as the optimal point of $\obj_k(u)$ on $\mathcal{U}$, $u_k^*$ satisfies $\big(u_k^* - \eta \nabla \obj_k(u_k^*)\big) - u_k^* \in \partial \psi_{\mathcal{U}}(u_k^*)$, 
    where $\partial \psi_{\mathcal{U}}(u_k^*)$ denotes the subdifferential of the indicator function $\psi_{\mathcal{U}}$ of $\mathcal{U}$ at $u_k^*$; \rev{(s.2) follows from the nonexpansiveness property of the projection operator, i.e., $\|\proj(x) - \proj(y)\| \leq \|x-y\|, \forall x,y \in \mathbb{R}^{p}$}; (s.3) uses the triangle inequality; (s.4) utilizes the $c_k$-Lipschitz continuity of the mapping $I-\eta\nabla \obj_k(\cdot)$ when $\obj_k$ is $\mu_k$-strongly convex and $L_k$-smooth.
    We plug the upper bound \eqref{eq:upper_bd_track_err_first_term} into \eqref{eq:upper_bd_track_err}, take expectations of both sides, and obtain
    \begin{align}\label{eq:track_err_recursive}
        \E_{v_{[T]}}&[\|w_{k+1} - u_{k+1}^*\|] \notag \\
            &\leq c_k\E_{v_{[T]}}[\|w_k - u_k^*\|] + \eta\E_{v_{[T]}}[\|\hat{\phi}_k - \nabla\obj_k(w_k)\|] \notag \\
            &\quad + \|u_{k+1}^* - u_k^*\|. 
    \end{align}
    By incorporating the upper bound \eqref{eq:upper_bd_dist_grad_overall} in \cref{lem:acc_hybrid_grad} and recursively applying \eqref{eq:track_err_recursive} for $k=0,\ldots,T-1$, we have
    \begin{align*}
        &\E_{v_{[T]}}[\|w_{T} - u_{T}^*\|] \\
            &\leq \Big(\!\prod_{k=0}^{T\!-\!1}c_k\!\Big) \cdot \|w_0 \!-\! u_0^*\| \!+\! \sum_{k\!=\!0}^{T\!-\!1}\! \Big(\!\prod_{s=k\!+\!1}^{T\!-\!1}c_s\!\Big) \cdot (\eta \gamma_k \!+\! \|u_{k+1}^* \!-\! u_k^*\|),
    \end{align*}
    where $\prod_{s=T}^{T-1} c_s$ is defined to be $1$. 
    Based on the parametric conditions, we know that $\forall k \in \mathbb{N}, c_k \in (0,1)$ and that $\rho = \max_k c_k \in (0,1)$. Moreover, $\sum_{k=0}^{T-1} \rho^{T-1-k} < 1/(1-\rho)$, and $\forall k \in \mathbb{N}, (1-\alpha_k) \in [0,1]$. Therefore, \eqref{eq:track_err_iss_upper_bd} holds.

\else
\fi

\subsection{Details of Numerical Experiments}\label{appendix:experiment_details}
    We describe the detailed setups of the numerical experiments in \cref{sec:experiment}. Our code is available \cite{he2024grayboxcode}. 

    \emph{\underline{Nonlinear dynamics \eqref{eq:sys_simulation}}} \par
    We draw the elements of the system matrices $A, B_1, B_2, C, D$ and $E$ in \eqref{eq:sys_simulation} from the normal distribution. We further scale $A$ so that its spectral radius is $0.05$, i.e., the dynamics are quickly contracting. The disturbances $d_x,d_y$ are sampled from the multivariate normal distribution, and their exact values are unknown beforehand.

    \emph{\underline{Static unconstrained problem \eqref{eq:opt_problem_simulation} and related controllers}} \par
    In the objective function of \eqref{eq:opt_problem_simulation}, $M_1 = M_3^{\top}M_3 \in \mathbb{R}^{15 \times 15}$ is positive definite, $m_2 \in \mathbb{R}^{15}$, $\lambda = 5\times 10^{-2}$, and the elements of $m_2$ and $M_3 \in \mathbb{R}^{15 \times 15}$ are drawn from normal distributions.

    The approximate sensitivity $\hat{H}$ used by various controllers is generated by perturbing $\left[C(I-A)^{-1}B_1\right]^\top$ with uniform noise, where the noise bound is $15\%$ of the average magnitude of the elements of the unperturbed matrix. That is, $\hat{H}$ is a noisy estimate of the sensitivity corresponding to the linear part of \eqref{eq:sys_simulation}.

    The step sizes for the model-based controllers with $H_k$, $\hat{H}$, and sensitivity learning\cite{picallo2021adaptive} as well as the stochastic extremum seeking algorithm\cite{liu2016stochastic} are $\eta = 5 \times 10^{-4}$. For the model-free and gray-box controllers, we set the step size as $10^{-4}$ and $2.5\times 10^{-4}$, respectively. The smoothing parameter used by both controllers is $\delta = 3 \times 10^{-3}$. The gray-box controller uses the rule \eqref{eq:adaptive_comb_coeff_case_constant} to determine the combination coefficients, and the constant $C$ therein is $5$.
    
    \emph{\underline{Time-varying constrained problem \eqref{eq:opt_problem_simulation_TV} and related controllers}} \par
    We generate the bounds $\underline{u}$ and $\bar{u}$ on the input $u$ from the multivariate normal distribution. Problem~\eqref{eq:opt_problem_simulation_TV} is time-varying in that every $10^3$ iterations, the positive definite $M_{1,k}$ and the coefficient vector $m_{2,k}$ in the objective $\Phi_k$ are regenerated from normal distributions, and the disturbances $d_{x,k},d_{y,k}$ are regenerated from uniform distributions.

    The approximate sensitivity $\hat{H}$ is a perturbed version of the sensitivity matrix $\left[C(I-A)^{-1}B_1 \right]^\top$. The perturbation noise follows uniform distributions.  

    The step sizes for the model-based controllers with $H_k$, $\hat{H}$, and sensitivity learning\cite{picallo2021adaptive} as well as the gray-box controller are $\eta = 5 \times 10^{-4}$. The model-free controller and the stochastic extremum seeking algorithm use step sizes of $10^{-4}$ and $2.5 \times 10^{-5}$, respectively. The smoothing parameter is $\delta = 0.05$. \revise{We set $C = 1.5$ in the rule \eqref{eq:bound_sensitivity_constant_TV} by which the gray-box controller tunes $\alpha_k$.}

    We obtain the comparator sequence $\{u_k^*\}$ (i.e., the optimal solutions to the time-varying problem~\eqref{eq:opt_problem_simulation_TV}) by calling the \textsf{fmincon} function of \textsf{MATLAB R2023b}.

    \else
        \balance
        \bibliographystyle{IEEEtran}
        \bibliography{article-abbrev}

@Article{chen2025continuous,
  author    = {Chen, Xin and others},
  journal   = {IEEE Trans. Autom. Control},
  title     = {Continuous-time zeroth-order dynamics with projection maps: Model-free feedback optimization with safety guarantees},
  year      = {2025},
  number    = {8},
  pages     = {5005-5020},
  volume    = {70},
  fjournal  = {IEEE Transactions on Automatic Control},
  publisher = {IEEE},
}

@article{scheinker2024100,
  title={100 years of extremum seeking: A survey},
  author={Scheinker, Alexander},
  journal={Automatica},
  volume={161},
  note={{A}rt.~no.~111481},
  year={2024},
  publisher={Elsevier}
}

@Article{simonetto2020time,
  author    = {Simonetto, Andrea and others},
  journal   = {Proc. IEEE},
  title     = {Time-varying convex optimization: Time-structured algorithms and applications},
  year      = {2020},
  number    = {11},
  pages     = {2032--2048},
  volume    = {108},
  fjournal  = {Proceedings of the IEEE},
  publisher = {IEEE},
}

@Article{simpson2021analysis,
  author    = {Simpson-Porco, John W},
  journal   = {IEEE Trans. Autom. Control},
  title     = {Analysis and synthesis of low-gain integral controllers for nonlinear systems},
  year      = {2021},
  number    = {9},
  pages     = {4148--4159},
  volume    = {66},
  fjournal  = {IEEE Transactions on Automatic Control},
  publisher = {IEEE},
}

@Article{bianchin2021time,
  author    = {Bianchin, Gianluca and others},
  journal   = {IEEE Trans. Control Netw. Syst.},
  title     = {Time-varying optimization of {LTI} systems via projected primal-dual gradient flows},
  year      = {2021},
  number    = {1},
  pages     = {474--486},
  volume    = {9},
  fjournal  = {IEEE Transactions on Control of Network Systems},
  publisher = {IEEE},
}

@article{cothren2022online,
  author={Cothren, Liliaokeawawa and others},
  //fjournal={IEEE Open Journal of Control Systems}, 
  journal={IEEE Open J. Control Syst.},
  title={Online Optimization of Dynamical Systems With Deep Learning Perception}, 
  year={2022},
  //volume={1},
  number={},
  //pages={306-321}
}

@article{simonetto2021personalized,
  author={Simonetto, Andrea and others},
  title={Personalized optimization with user's feedback},
  journal={Automatica},
  volume={131},
  //pages={109767},
  note = {{A}rt.~no.~109767},
  year={2021},
  publisher={Elsevier}
}

@article{belgioioso2022online,
  author={Belgioioso, Giuseppe and others},
  title={Online feedback equilibrium seeking},
  journal  = {IEEE Trans. Autom. Control},
  year={2025},
  volume={70},
  number={1},
  pages={203-218},
  fjournal = {IEEE Transactions on Automatic Control}
}

@Article{bianchin2021online,
  author   = {Bianchin, Gianluca and others},
  journal  = {IEEE Trans. Autom. Control},
  title    = {Online Stochastic Optimization for Unknown Linear Systems: Data-Driven Controller Synthesis and Analysis},
  year     = {2024},
  volume   = {69},
  number   = {7},
  pages    = {4411-4426},
  fjournal = {IEEE Transactions on Automatic Control},
}

@article{nonhoff2022online,
  title={Online convex optimization for data-driven control of dynamical systems},
  author={Nonhoff, Marko and others},
  //journal={IEEE Open Journal of Control Systems},
  journal={IEEE Open J. Control Syst.},
  //volume={1},
  //pages={180--193},
  year={2022},
  publisher={IEEE}
}

@Article{picallo2021adaptive,
  author    = {Picallo, Miguel and others},
  journal   = {Electr. Pow. Syst. Res.},
  title     = {Adaptive real-time grid operation via Online Feedback Optimization with sensitivity estimation},
  year      = {2022},
  note      = {{A}rt.~no.~108405},
  volume    = {212},
  //number  = {108405},
  fjournal  = {Electric Power Systems Research},
  publisher = {Elsevier},
}

@INPROCEEDINGS{dominguez2023online,
  title={An Online Feedback Optimization Approach to Voltage Regulation in Inverter-Based Power Distribution Networks},
  author={Dominguez-Garcia, Alejandro D and others},
  booktitle = {Proc. Amer. Control Conf.},
  //address = {San Diego, CA, USA},
  pages={1868-1873},
  year={2023}
}

@Article{krishnamoorthy2023model,
  author    = {Krishnamoorthy, Dinesh and others},
  journal   = {AlChE J.},
  title     = {Model-free real-time optimization of process systems using safe {B}ayesian optimization},
  year      = {2023},
  number    = {4},
  //pages     = {e17993},
  note = {{A}rt.~no.~e17993},
  volume    = {69},
  fjournal  = {AIChE Journal},
  publisher = {Wiley Online Library},
}

@inproceedings{chen2020model,
  author={Chen, Yue and others},
  title={Model-free primal-dual methods for network optimization with application to real-time optimal power flow},
  booktitle={Proc. Amer. Control Conf.},
  pages={3140--3147},
  year={2020}
}

@article{tang2023zeroth,
  title={Zeroth-order feedback optimization for cooperative multi-agent systems},
  author={Tang, Yujie and others},
  journal={Automatica},
  volume={148},
  //pages={110741},
  note = {{A}rt.~no.~110741},
  year={2023},
  publisher={Elsevier}
}

@Article{poveda2017robust,
  author    = {Poveda, Jorge I and others},
  journal   = {IEEE Trans. Autom. Control},
  title     = {A robust event-triggered approach for fast sampled-data extremization and learning},
  year      = {2017},
  number    = {10},
  pages     = {4949--4964},
  volume    = {62},
  fjournal  = {IEEE Transactions on Automatic Control},
  publisher = {IEEE},
}

@inproceedings{qu2021exploiting,
  author={Qu, Guannan and others},
  title={Exploiting linear models for model-free nonlinear control: A provably convergent policy gradient approach},
  booktitle={Proc. IEEE 60th Conf. Decis. Control},
  pages={6539--6546},
  year={2021}
}

@Article{ma2023reinforcement,
  author    = {Ma, Hao and others},
  journal   = {Auton. Robot.},
  title     = {Reinforcement learning with model-based feedforward inputs for robotic table tennis},
  year      = {2023},
  number    = {8},
  pages     = {1387--1403},
  volume    = {47},
  fjournal  = {Autonomous Robots},
  publisher = {Springer},
}

@article{li2023certifying,
  author={Li, Tongxin and others},
  title={Certifying Black-Box Policies With Stability for Nonlinear Control},
  journal={IEEE Open J. Control Syst.},
  volume={2},
  pages={49--62},
  year={2023},
  publisher={IEEE}
}

@book{dontchev2009implicit,
  title={Implicit functions and solution mappings},
  author={Dontchev, Asen L and others},
  volume={543},
  year={2009},
  publisher={Springer},
  address={New York, NY, USA}
}

@Article{gao2018information,
  author    = {Gao, Xiang and others},
  journal   = {J. Sci. Comput.},
  title     = {On the information-adaptive variants of the {ADMM}: an iteration complexity perspective},
  year      = {2018},
  number    = {1},
  pages     = {327--363},
  volume    = {76},
  fjournal  = {Journal of Scientific Computing},
  publisher = {Springer},
}

@Article{nesterov2017random,
  author    = {Nesterov, Yurii and others},
  journal   = {Found. Comput. Math.},
  title     = {Random gradient-free minimization of convex functions},
  year      = {2017},
  number    = {2},
  pages     = {527--566},
  volume    = {17},
  fjournal  = {Foundations of Computational Mathematics},
  publisher = {Springer},
}

@Article{he2022model,
  author   = {He, Zhiyu and others},
  journal  = {IEEE Trans. Autom. Control},
  title    = {Model-Free Nonlinear Feedback Optimization},
  year     = {2024},
  volume   = {69},
  number   = {7},
  pages    = {4554-4569},
  fjournal = {IEEE Transactions on Automatic Control},
}

@Article{ospina2022feedback,
  author    = {Ospina, Ana M and others},
  journal   = {IEEE Control Syst. Lett.},
  title     = {Feedback-based optimization with sub-weibull gradient errors and intermittent updates},
  year      = {2022},
  pages     = {2521--2526},
  volume    = {6},
  fjournal  = {IEEE Control Systems Letters},
  publisher = {IEEE},
}

@article{zhang2022new,
  title={A new one-point residual-feedback oracle for black-box learning and control},
  author={Zhang, Yan and others},
  journal={Automatica},
  volume={136},
  //pages={110006},
  note = {{A}rt.~no.~110006},
  year={2022},
  publisher={Elsevier}
}

@Article{zhao2021bandit,
  author   = {Zhao, Peng and others},
  journal  = {J. Mach. Learn. Res.},
  title    = {Bandit convex optimization in non-stationary environments},
  year     = {2021},
  number   = {1},
  pages    = {5562--5606},
  volume   = {22},
  fjournal = {Journal of Machine Learning Research},
}

@book{hazan2022introduction,
  title={Introduction to online convex optimization},
  author={Hazan, Elad},
  year={2022},
  publisher={MIT Press},
  address = {Princeton, NJ}
}

@inproceedings{hazan2017efficient,
  title={Efficient regret minimization in non-convex games},
  author={Hazan, Elad and others},
  booktitle={Proc. Int. Conf. Mach. Learn.},
  //booktitle={International Conference on Machine Learning},
  pages={1433--1441},
  year={2017}
}

@Article{ghadimi2016mini,
  author    = {Ghadimi, Saeed and others},
  journal   = {Math. Program.},
  title     = {Mini-batch stochastic approximation methods for nonconvex stochastic composite optimization},
  year      = {2016},
  number    = {1},
  pages     = {267--305},
  volume    = {155},
  fjournal  = {Mathematical Programming},
  publisher = {Springer},
}

@Article{liu2016stochastic,
  author    = {Liu, Shu-Jun and others},
  journal   = {IEEE Trans. Autom. Control},
  title     = {Stochastic averaging in discrete time and its applications to extremum seeking},
  year      = {2016},
  number    = {1},
  pages     = {90--102},
  volume    = {61},
  fjournal  = {IEEE Transactions on Automatic Control},
  publisher = {IEEE},
}

@article{he2024gray,
  author={He, Zhiyu and others},
  title={Gray-box nonlinear feedback optimization},
  journal={arXiv preprint arXiv:2404.04355},
  year={2024}
}

@misc{he2024grayboxcode,
  author       = {He, Zhiyu},
  title        = {Gray-Box Nonlinear Feedback Optimization},
  year         = {2024},
  howpublished = {\textit{GitHub Repository}},
  url          = {https://github.com/zyhe/Gray-box-nonlinear-feedback-optimization},
  note         = {Accessed: Jun. 15, 2026}
}
    \fi
\end{document}